\documentclass[letterpaper, 11pt, reqno]{amsart}

\usepackage{fullpage}
\usepackage{graphicx}
\usepackage{amsmath, amssymb, amsfonts, amsthm, mathtools}
\usepackage{url}
\usepackage[foot]{amsaddr}

\usepackage{color, comment}
\usepackage{xcolor} 
\definecolor{darkblue}{rgb}{0,0,0.4}
\usepackage[plainpages=false, colorlinks=true, citecolor=blue, filecolor=black, linkcolor=red, urlcolor=darkblue]{hyperref}

\newtheorem{thm}{Theorem}[section]
\newtheorem{prop}[thm]{Proposition}
\newtheorem{conj}[thm]{Conjecture}
\newtheorem{lem}[thm]{Lemma}
\theoremstyle{remark}
\newtheorem{rem}[thm]{Remark}
\newtheorem{defin}[thm]{Definition}

\newcommand{\calA}{\mathcal{A} }
\newcommand{\calI}{\mathcal{I} }
\newcommand{\ws}{\overset{\ast}{\rightharpoonup}}
\newcommand{\del}{\partial} 
\newcommand{\R}{\mathbb{R}} 

\newcommand{\N}{\mathbb{N}} 
\newcommand{\D}{\mathbb{D}} 
\renewcommand{\S}{\mathbb{S}} 
\newcommand{\vareps}{\varepsilon} 
\newcommand{\abs}[1]{\left\lvert #1 \right\rvert} 

\DeclareMathOperator{\tr}{tr}

\usepackage{dsfont}
\newcommand{\one}{\mathds{1}}

\graphicspath{{./figs/}} 

\usepackage[backend=biber, giveninits=true, style=numeric, natbib=true, url=false, sorting=nyt, doi=true, backref=false]{biblatex}
\renewbibmacro{in:}{\ifentrytype{article}{}{\printtext{\bibstring{in}\intitlepunct}}}

\AtEveryBibitem{\clearfield{month}}
\AtEveryBibitem{\clearfield{day}}
\AtEveryBibitem{\clearfield{pages}}
\AtEveryBibitem{\clearfield{volume}}
\AtEveryBibitem{\clearfield{number}}

\begin{document}
\title{Extremal Steklov-Neumann eigenvalues}
\thanks{C.-Y.K. acknowledges support from NSF DMS-2208373 and DMS-2513176. B.O. acknowledges support from NSF DMS-2136198 and DMS-2513175. C.H.T. acknowledges support from the National Science and Technology Council of Taiwan (NSTC 114-2115-M-110-001-MY2).}

\author{Chiu-Yen Kao}
\address{Department of Mathematical Sciences, Claremont McKenna College, Claremont, CA 91711, USA}
\email{ckao@cmc.edu}

\author{Braxton Osting}
\address{Department of Mathematics, University of Utah, Salt Lake City, UT 84112, USA}
\email{osting@math.utah.edu}

\author{Chee Han Tan}
\address{Department of Applied Mathematics, National Sun Yat-sen University, Kaohsiung 80424, Taiwan}
\email{tanch@math.nsysu.edu.tw}

\author{Robert Viator}
\address{Department of Mathematics, Denison University, Granville, OH 43023, USA}
\email{viatorr@denison.edu}

\subjclass[2020]{
31A25,  
35P15, 
49M41, 
65K10, 
65N25, 
49R05, 
}

\keywords{Steklov-Neumann eigenproblem, weighted Steklov eigenproblem, extremal eigenvalues, variational analysis, optimal design}

\date{\today}

\begin{abstract} 
Let $\Omega$ be a bounded open planar domain with smooth connected boundary, $\Gamma$, that has been partitioned into two disjoint components,  $\Gamma = \Gamma_S \sqcup \Gamma_N$. We  consider the Steklov-Neumann eigenproblem on $\Omega$, where a harmonic function is sought that satisfies the Steklov boundary condition on $\Gamma_S$ and the Neumann boundary condition on $\Gamma_N$. We pose the extremal eigenvalue problems (EEPs) of minimizing/maximizing the $k$-th non-trivial Steklov-Neumann eigenvalue among boundary partitions of prescribed measure. We formulate a relaxation of these EEPs in terms of weighted Steklov eigenvalues where an $L^\infty(\Gamma)$ density replaces the boundary partition. For these relaxed EEPs, we establish existence and prove optimality conditions. We also prove a homogenization result that allows us to use solutions to the relaxed EEPs to infer properties of solutions to the original EEPs. For a disk, we provide numerical and asymptotic evidence that the minimizing arrangement of $\Gamma_S\sqcup \Gamma_N$ for the $k$-th eigenvalue consists of $k+1$ connected components that are symmetrically arranged on the boundary. For a disk, for $k = 1$, the constant density is a maximizer for the relaxed problem; we also provide numerical and asymptotic evidence that for $k\ge 2$, the maximizing density for the relaxed problem is a non-trivial function; a sequence of rapidly oscillating Steklov/Neumann boundary conditions approach the supremum value. 
\end{abstract}

\maketitle

\section{Introduction} 
Fix a bounded open domain $\Omega\subset\R^2$ with 
Lipschitz connected boundary $\del\Omega = \Gamma$. For a disjoint partition of the boundary, $\Gamma =   \Gamma_S \sqcup \Gamma_N$, we consider the Steklov-Neumann (S-N) eigenproblem: 
\begin{subequations} \label{e:Steklov}
\begin{alignat}{2}
\Delta u & = 0 && \qquad \text{ in $\Omega$}, \\
\del_n u & = \sigma u && \qquad \text{ on $\Gamma_S$}, \\
\del_n u & = 0 && \qquad \text{ on $\Gamma_N$}. 
\end{alignat}
\end{subequations}
Here, $\Delta =  \partial^2_{x}+\partial^2_{y}$ is the Laplacian and $\partial_n$ is the normal derivative on the boundary. 
The S-N eigenvalues are countable and we enumerate them, counting multiplicity, in increasing order, 
$ 0 = \sigma_0 < \sigma_1 \leq \sigma_2\leq \dots \nearrow\infty$. 
The S-N eigenvalues satisfy a variational principle,
$$
\sigma_k = \min_{E_k \subset H^1(\Omega)} \ \max_{u \in E_k\setminus\{0\}} \  \frac{\int_\Omega \abs{\nabla u}^2 dx }{ \int_{\Gamma_S} u^2 \, ds }, 
$$ 
where $E_k$ is a $(k + 1)$--dimensional subspace of the Sobolev space,  $H^1(\Omega)$.

For a given $\alpha \in (0,1]$ and $k \in \N = \{1,2,\ldots\}$, we consider the S-N eigenvalue optimization problem of finding the disjoint partition of the boundary, $\Gamma = \Gamma_S \sqcup \Gamma_N$ with $\abs{\Gamma_S} = \alpha \abs{\Gamma}$, that minimizes the $k$-th non-trivial eigenvalue,  
\begin{equation} \label{e:Min}
\sigma_{k,\alpha}^{\bigtriangledown} = \inf_{\substack{\Gamma_S\subset \Gamma \\ \abs{\Gamma_S} = \alpha\abs{\Gamma}}} \ \sigma_k(\Gamma_S).
\end{equation}
We also consider maximizing the $k$-th non-trivial eigenvalue, 
\begin{equation} \label{e:Max}
\sigma_{k,\alpha}^{\triangle} = \sup_{\substack{\Gamma_S\subset \Gamma \\ \abs{\Gamma_S} = \alpha\abs{\Gamma}}} \ \sigma_k(\Gamma_S).
\end{equation}

The optimization problems \eqref{e:Min} and \eqref{e:Max} were posed by Nilima Nigam in the AIM Problem List (AimPL) on Steklov eigenproblems \cite[Problem 4.1]{AimPL}. Steklov eigenvalues appear in many applications \cite{Kuznetsov_2014, Girouard_2017}, so the extremal eigenvalue problems  \eqref{e:Min} and \eqref{e:Max} are of fundamental interest. We offer a physical interpretation in terms of the analogous extremal Laplacian eigenvalue problems with mixed Dirichlet/Neumann boundary conditions (B.C.) studied by  S.~J.~Cox and P.~X.~Uhlig \cite{Cox_2003}. These problems, for $k=1$, have the interpretation that if we are allowed to fasten a prescribed fraction of the boundary of a drum (corresponding to Dirichlet B.C.), which portion should we fasten to produce the lowest/highest bass note (first Laplacian eigenvalue)? The S-N eigenvalue problem \eqref{e:Steklov} is the limit of this eigenvalue problem where the mass of the drum is concentrated at $\Gamma_S$ \cite{Lamberti_2015}. Thus, the optimization problem \eqref{e:Min} for $k = 1$ has the following physical interpretation: Suppose we are allowed to fasten a fixed fraction of the boundary of a drum (corresponding to Steklov B.C.), where the mass of the drum will also be concentrated. Which portion should we fasten to produce the lowest bass note? 

The S-N eigenvalue problem has attracted the interest of researchers in recent years. In \cite{Banuelos2010}, Ba{\~n}uelos, Kulczycki, Polterovich, and Siudeja established inequalities between mixed S-N and mixed Steklov-Dirichlet eigenvalues for certain domains satisfying the weak John's condition. These inequalities are used to obtain geometric information about the nodal sets of S-N eigenfunctions.  
In \cite {Ammari_2020}, Ammari, Imeri, and Nigam analyzed how small changes in the S-N partition move the spectrum and used the resulting shape-derivative formulas to design a boundary-integral method that tunes mixed S-N conditions to hit target resonances. They also treated multiple eigenvalues and validated the asymptotics numerically on canonical shapes. 
In \cite{Arias_2024}, Arias-Marco, Dryden, Gordon, Hassannezhad, Ray, and Stanhope proved a sharp isoperimetric inequality for mixed Steklov eigenvalues on surfaces by relating S-N and Steklov–Dirichlet spectra, developed full asymptotics for mixed problems under mild interface regularity, and studied symmetry-constrained maximization, showing that the upper bounds converge to the Hersch-Payne-Schiffer asymptotics as $k \to \infty$. 
In \cite{Grebenkov_2025}, 
Grebenkov gave a ``small-window'' asymptotic analysis for the mixed S-N problem (framed via an integral-operator reformulation). The asymptotics are checked numerically for an arc on a disk and a spherical cap on a ball, and are linked to diffusion-controlled reactions (mean first-reaction times and reactivity-dependent corrections).
In \cite{Basak_2024}, Basak and Verma studied Steklov and mixed S-N spectra on domains with a spherical hole. They obtained an explicit expression for the second nonzero Steklov eigenvalue on an annulus and derived sharp upper bounds for the first $n$ nonzero Steklov eigenvalues on 4--fold symmetric domains in $\R^n$, with analogous bounds for the mixed problem; examples show the symmetry hypothesis is essential.

\subsection*{Weighted Steklov eigenproblem and relaxed extremal S-N eigenvalues}
For $\rho \in L^\infty(\Gamma)$ and $\rho(s) \geq 0$ for a.e. $s\in \Gamma$, we ``relax'' the S-N eigenproblem
 \eqref{e:Steklov} and consider the \emph{weighted Steklov eigenproblem}, 
\begin{subequations} \label{e:relax}
\begin{alignat}{2}
\label{e:relaxa} \Delta u & = 0 && \qquad \text{ in $\Omega$}, \\
\label{e:relaxb} \del_n u & = \lambda\rho u && \qquad \text{ on $\Gamma$}. 
\end{alignat}
\end{subequations}
Observe that we recover the SN eigenproblem \eqref{e:Steklov} for $\rho = \one_{\Gamma_S}$, the indicator function on $\Gamma_S$. As for the S-N eigenproblem, the weighted Steklov eigenvalues are countable and we enumerate them, counting multiplicity, in increasing order, 
$ 0 = \lambda_0 < \lambda_1 \leq \lambda_2\leq \dots \nearrow\infty$.  
The weighted Steklov eigenvalues satisfy a variational principle, 
\begin{equation}
\label{e:VarPrinc}
\lambda_k = \min_{E_k \subset H^1(\Omega)} \ \max_{u \in E_k\setminus\{0\}}  \ \frac{\int_\Omega \abs{\nabla u}^2 dx }{\int_\Gamma \rho u^2\, ds}. 
\end{equation}
This weighted Steklov eigenproblem is of interest in its own right, as a conformal mapping applied to the Steklov eigenvalue problem will introduce a boundary weight \cite{weinstock1954inequalities, alhejaili2019numerical, Oudet_2021}. 

For fixed $\alpha \in (0,1]$, define the admissible set of densities, $\calA_\alpha$, by 
\begin{equation} \label{e:Ad}
\calA_\alpha\coloneqq \left\{\rho \in L^\infty(\Gamma) \colon 0 \leq \rho(s) \leq 1 \textrm{ for a.e. } s \in \Gamma, \, \int_{\Gamma} \rho(s)\, ds = \alpha \abs{\Gamma} \right\}. 
\end{equation}
For fixed $\alpha \in (0,1]$ and $k\in\N$, we then consider relaxations of the minimization problem \eqref{e:Min}, 
\begin{equation} \label{e:RelaxMin}
\lambda_{k,\alpha}^{\bigtriangledown} = \inf_{\rho\in\calA_\alpha} \ \lambda_k(\rho), 
\end{equation}
and the maximization problem \eqref{e:Max}, 
\begin{equation} \label{e:RelaxMax}
\lambda_{k,\alpha}^{\triangle} = \sup_{\rho\in\calA_\alpha} \ \lambda_k(\rho). 
\end{equation}
Clearly, we have that
$
\lambda_{k,\alpha}^{\bigtriangledown} 
\leq \sigma_{k,\alpha}^{\bigtriangledown}
\leq \sigma_{k,\alpha}^{\triangle} 
\leq \lambda_{k,\alpha}^{\triangle} $. Note that these are all finite, since  by the Hersch-Payne-Schiffer inequality \cite{Hersch_1974,Girouard_2010}
\begin{equation}
\label{e:UpperBnd}
\lambda_k(\rho) \leq \frac{2\pi}{\alpha\abs{\Gamma}} k, \qquad \forall \rho \in \calA_\alpha.  
\end{equation}

\subsection*{Main results}
We first establish the following existence results. 

\begin{thm}(Existence)  
\label{p:existence}
For every $k\in \N$ and $\alpha \in (0,1]$, there exists $\rho_k^{\bigtriangledown} \in \calA_\alpha$ that attains the infimum in \eqref{e:RelaxMin} and there exists $\rho_k^{\triangle} \in \calA_\alpha$ that attains the supremum in \eqref{e:RelaxMax}. 
\end{thm}
Theorem~\ref{p:existence} is proven in Section~\ref{s:Proof}. The statement follows from the continuity of $\rho\mapsto\lambda_k(\rho)$ on $\calA_\alpha$ and the compactness of $\calA_\alpha$ with respect to the weak$^*$ topology on $L^\infty(\Gamma)$.

\begin{rem}
\label{r:nonconvex}
The map $\rho \mapsto \lambda_k^{-1}(\rho)$ is not convex for all $k\ge 1$; see \cite[Figure~1]{KaoAbbas2025}. This seems to be misunderstood in the literature for eigenvalues of other elliptic operators. In Appendix~\ref{a:counterexample}, we give a simple example of a Sturm-Liouville operator that has non-convex higher eigenvalues.
\end{rem} 

Our next result states that weighted Steklov eigenvalues for a rapidly oscillating characteristic microstructure converge to those of the homogenized density. 
\begin{thm}(Homogenization)  
\label{p:homogenization} 
Suppose $\rho\in L^\infty(\Gamma)$ and $\rho_n = \one_{S_n}$ for measurable $S_n \subset \Gamma$, $n\geq 1$ such that $\lim\limits_{n\to \infty} \abs{A\cap S_n}  = \int_A \rho(s) \, ds$ for every measurable set $A \subset \Gamma$. Then $\rho_n\ws\rho$ in $L^\infty(\Gamma)$ and for fixed $k\geq 0$, $\lambda_k(\rho_n) \to \lambda_k(\rho)$.   
Consequently, $\lambda_{k,\alpha}^{\bigtriangledown} = \sigma_{k,\alpha}^{\bigtriangledown}$ and $\sigma_{k,\alpha}^{\triangle} = \lambda_{k,\alpha}^{\triangle}$. 
\end{thm}
Theorem~\ref{p:homogenization} is proven in Section~\ref{s:Proof} and has the following consequence. The inf/sup in \eqref{e:Min}/\eqref{e:Max} is attained by an admissible boundary partition $\Gamma = \Gamma_S \sqcup \Gamma_N$ if and only if there is a bang-bang extremal density for \eqref{e:RelaxMin}/\eqref{e:RelaxMax}. Otherwise, there is a sequence of rapidly oscillating characteristic microstructures, $\rho_n$, that weak* converge to the extremal density, $\rho^*$, such that $\lambda_k(\rho_n) \to \lambda_k(\rho^*)$. Here, we say that $\rho \in \calA_\alpha$ is a ``bang-bang'' density if $\rho = \one_S$ for a measurable set $S \subseteq \Gamma$.

Theorem~\ref{p:existence} and the variational principle \eqref{e:VarPrinc} for $\lambda_k$ motivate the following heuristics: 
\begin{itemize}
\item[(i)] to minimize $\lambda_k$, we should put $\rho= 1$ (i.e., $\Gamma_S$) where $\abs{u_k}$ is large, and
\item[(ii)] to maximize $\lambda_k$, we should put $\rho= 1$ (i.e.,  $\Gamma_S$) where $\abs{u_k}$ is small. 
\end{itemize} 
Of course, this heuristic only makes sense when the reference to $u_k$ is well-defined, (i.e.,  when $\lambda_k$ is simple). The following theorem generalizes this intuition to a necessary optimality condition for a (possibly) higher multiplicity eigenvalue. We introduce the notation $[n]\coloneqq\{1, 2, \ldots, n\}$ for $n\in\N$ and define $L^2(\Gamma; \rho)$ to be the weighted $L^2$ space on $\Gamma$ with a nonnegative weight $\rho$. 

\begin{thm}(Optimality) 
\label{prop:OptCondMultEig}
Fix $k \in\N$ and $\alpha \in (0,1]$. 
Let $\rho_k^\bigtriangledown\in \calA_\alpha$ be a local minimizer for \eqref{e:RelaxMin}. 
Assume $\lambda_k(\rho_k^\bigtriangledown)$ is a multiplicity $m_k$ eigenvalue and let $\{u_{i, k}\}_{i \in [m_k]}$ be an $L^2(\Gamma; \rho_k^\bigtriangledown)$-orthonormal choice of eigenfunctions. Then there exists $c > 0$ such that the squared eigenfunction boundary values satisfy 
\begin{equation} 
\label{e:minRhoCondMult} 
\sum_{i\in [m_k]} u_{i,k}^2\Big|_{\Gamma}(s) 
\begin{cases} 
\, \leq c, & \textrm{if } \rho_k^\bigtriangledown(s) = 0, \\ 
\, = c, & \textrm{if } \rho_k^\bigtriangledown(s) \in (0, 1), \\  
\, \geq c, &\textrm{if } \rho_k^\bigtriangledown(s) = 1.
\end{cases} 
\end{equation}

Similarly, let $\rho_k^\triangle\in\calA_\alpha$ be a local maximizer for \eqref{e:RelaxMax}. 
Assume $\lambda_k(\rho_k^\triangle)$ is a multiplicity $m_k$ eigenvalue and let $\{u_{i, k}\}_{i \in [m_k]}$ be an $L^2(\Gamma; \rho_k^\triangle)$-orthonormal choice of eigenfunctions.
Then there exists $c > 0$ such that the squared eigenfunction boundary values satisfy 
\begin{equation} 
\label{e:maxRhoCondMult} 
\sum_{i\in [m_k]} u_{i,k}^2\Big|_{\Gamma}(s) 
\begin{cases} 
\, \geq c, & \textrm{if } \rho_k^\triangle(s) = 0, \\ 
\, = c, & \textrm{if } \rho_k^\triangle(s) \in (0, 1), \\  
\, \leq c, &\textrm{if } \rho_k^\triangle(s) = 1.
\end{cases} 
\end{equation}
\end{thm}

The constant $c>0$ in the optimality conditions \eqref{e:minRhoCondMult} and \eqref{e:maxRhoCondMult} depends on the choice of eigenbasis $\{ u_{i,k}\}_{i\in [m_k]}$; see the proof in Section \ref{s:Proof} for details.  Interestingly, for a simple extremal eigenvalue (either minimizer or maximizer), the corresponding eigenfunction is constant on any interval for which the extremal density takes values between 0 and 1; see Figures~\ref{f:EvenkMaximizer} and \ref{f:SmallkMaximizer}. Theorem~\ref{prop:OptCondMultEig} is proven in Section~\ref{s:Proof}. 
The proof uses a local variational argument relying on the Gateaux derivative of a non-simple eigenvalue; see Lemma~\ref{l:PropSigma}(e). 
This local analysis is similar to asymptotics of Steklov eigenvalues with respect to shape perturbations \cite{Viator_2018,Viator_2020,Viator_2022,Tan_2023,Schroeder_2023,Tan_2024,Alland_2026}.

In Section~\ref{s:disk}, we consider the case when $\Omega$ is a unit disk $\D\subset\R^2$, so that $\Gamma = \partial \Omega = \S^1$. 
\begin{defin} \label{d:Sn}
Let $\calI_\alpha^{n}\subset \S^1$ be the union of $n$ evenly-spaced equal-length intervals, given by 
\begin{equation*} 
\calI_\alpha^{n}= \bigcup_{\ell \in [n]} \calI_\ell, \quad \calI_\ell = \left[ \frac{\pi (2 \ell - \alpha)}{n}, \ \frac{\pi (2 \ell + \alpha)}{n}\right], \quad \ell \in [n]\coloneqq\{1, 2, \ldots,n\}. 
\end{equation*}
Each interval $\calI_\ell$ has length $\frac{2\pi \alpha}{n}$, so that  $\abs{\calI_\alpha^n} = 2\pi\alpha = \alpha\abs{\Gamma}$. The indicator function  $\one_{\calI_\alpha^{n}}$ satisfies $\int_0^{2\pi} \one_{\calI_\alpha^{n}}(\theta)\, d\theta = \alpha\abs{\Gamma}$, so is an admissible density; see \eqref{e:Ad}. 
\end{defin}

We first make the following conjecture. We will prove some of the statements in the conjecture, as stated in the subsequent theorem and discussed further with computational results in Section~\ref{s:CompRes}. 

\begin{conj}(Disk) \label{c:Disk}
Let $\Omega = \D$. We conjecture the following to hold. 
\begin{enumerate}
\item (Minimizers) For $k\geq 1$ and $\alpha \in (0,1]$, the unique (up to rotation) minimizer of $\lambda_k$ in \eqref{e:RelaxMin} is the indicator function $\rho_k^\bigtriangledown = \one_{\calI_\alpha^{k+1}}$, where $\calI_\alpha^{k + 1}\subset\S^1$ is given in Definition~\ref{d:Sn}. For odd $k$, the multiplicity is one and for even $k$, the multiplicity is two. Approximate numerical solutions are shown in Figures \ref{f:kOdd}, \ref{f:kEven}, and \ref{f: min optimiality}. \\
Consequently, for $k\geq 1$, the unique (up to rotation) minimizer of $\sigma_k$ in \eqref{e:Min} is $\calI_\alpha^{k+1} \subset \S^1$.  

\item (Maximizers) For $k = 1$  and $\alpha \in (0,1]$, the constant density $\rho_1^{\triangle} = \alpha$ is a maximizer of $\lambda_1$ in \eqref{e:RelaxMax}. 
For $k\geq 2$, the unique (up to rotation) maximizer of $\lambda_k$ in \eqref{e:RelaxMax} has the form in 
Figures~\ref{f:EvenkMaximizer} and \ref{f:SmallkMaximizer}; 
$\rho_k^{\triangle}$ is the $k$-repetition of an interval that is constant ($=1$) on a first subinterval and is a non-trivial, strictly positive, and convex function on a second subinterval.
For even $k$, the multiplicity is one (consequently, the eigenfunction is constant on the second subintervals where $\rho_k^{\triangle} < 1$) and for odd $k$, the multiplicity is two. \\ 
For $k\geq 1$, there are configurations consisting of rapidly oscillating Steklov/Neumann B.C. that can arbitrarily well approximate the supremum value in \eqref{e:Max}.  
\end{enumerate}
\end{conj} 

Note that despite the objective function in \eqref{e:RelaxMin} being neither convex nor concave (see Remark~\ref{r:nonconvex}), we conjecture the minimizers in \eqref{e:RelaxMin} to be bang-bang solutions; see Figures~\ref{f:kOdd}, \ref{f:kEven}, and \ref{f: min optimiality}. 
\begin{thm}(Disk) \label{p:Disk}
Let $\Omega = \D$. The following statements hold. 
\begin{enumerate}
    \item (Maximizer of $\lambda_1$) For $k = 1$, the constant density $\rho_1^{\triangle} = \alpha$ is a maximizer for \eqref{e:RelaxMax}. 

    \item (Local analysis for constant density) 
    For even $k\geq 2$  and $\alpha \in (0,1]$, the constant density $\rho = \alpha$ is a critical point for \eqref{e:RelaxMin}, in the sense that there is no admissible perturbation with negative Gateaux derivative. 
    For odd $k\ge 1$  and $\alpha \in (0,1]$, the constant density is not a local minimizer for \eqref{e:RelaxMin}.\\
    For odd $k\ge 1$  and $\alpha \in (0,1]$, the constant density is a critical point for \eqref{e:RelaxMax}, in the sense that there is no admissible perturbation with positive Gateaux derivative.  
    For even $k\ge 2$  and $\alpha \in (0,1]$, the constant density is not a local maximizer for \eqref{e:RelaxMax}.  
    
    \item (Asymptotic results for $\alpha \approx 1$) 
    Let $k \geq 1$  and $\vareps > 0$. Consider the extremal eigenvalue problems \eqref{e:Min} and \eqref{e:Max} for $\abs{\Gamma_S} = (1 - \vareps)\abs{\Gamma}$ with $\vareps\ll 1$, so that $\abs{\Gamma_N} = \vareps\abs{\Gamma}$. Let $\calI_{1 - \vareps}^{n}\subset\S^1$ be defined as in Definition~\ref{d:Sn}. 
    For fixed $k\ge 1$, in the limit as $\vareps\to 0$, $\Gamma_S = \calI_{1 - \vareps}^{k + 1}$ is a local minimizer for \eqref{e:Min} and $\Gamma_S = \calI_{1 - \vareps}^k$ is a local maximizer for \eqref{e:Max}.  
    \item (Homogenization) 
    Let $\calI_\alpha^{n}\subset\S^1$ be defined as in Definition~\ref{d:Sn}.
    The eigenvalues of the S-N eigenproblem \eqref{e:Steklov} with $\Gamma_S = \calI_\alpha^n$ converge, as $n \to \infty$, to the eigenvalues of the weighted Steklov eigenproblem \eqref{e:relax} with constant density $\rho = \alpha$.    
\end{enumerate}
\end{thm}
Theorem~\ref{p:Disk}(2)--(3) is proven in Section~\ref{s:disk}. 
Theorem~\ref{p:Disk}(1) is a special case of Weinstock's inequality \cite{weinstock1954inequalities, Girouard_2010}. 
Theorem~\ref{p:Disk}(2) rules out the possibility that the constant density is the minimizer/maximizer for \eqref{e:RelaxMin}/\eqref{e:RelaxMax} for $k\geq 2$. 
Theorem~\ref{p:Disk}(3) establishes Conjecture~\ref{c:Disk}(1) in the asymptotic limit when $\alpha \to 1$. This is stated somewhat informally; by ``in the limit as $\vareps\to 0$'' we mean that the linearization of the problems about $\vareps = 0$ have the stated properties. 
Theorem~\ref{p:Disk}(4) shows rapidly oscillating Steklov/Neumann B.C. homogenizes and behaves like a constant density and follows directly from Theorem~\ref{p:homogenization}.  
Based on Conjecture~\ref{c:Disk}(2) concerning $\rho_k^{\triangle}$ and Theorem~\ref{p:homogenization},
there is not a partition $\Gamma = \Gamma_S \sqcup \Gamma_N$ that attains the supremum \eqref{e:Max}  for $k\geq 1$.

In Section~\ref{s:CompMeth}, we describe a computational approach for solving \eqref{e:RelaxMin} and \eqref{e:RelaxMax}. 
We describe an efficient integral equation approach for solving the weighted Steklov eigenproblem \eqref{e:relax} based on representing the eigenfunction in terms of a modified single layer potential. This formulation uses the weak formulation of the integral eigenvalue equation to extend the method in \cite{Akhmetgaliyev2016} to include boundary weights. We also describe a gradient-based optimization method for approximating the solutions to the non-convex, non-smooth extremal eigenvalue problems  \eqref{e:RelaxMin} and \eqref{e:RelaxMax}. 
In Section~\ref{s:CompRes}, these computational methods are used in several computational experiments.
For the disk $\Omega = \D$, we observe that the extremal densities are very structured, as summarized in Conjecture~\ref{c:Disk}. 

\subsection*{Outline} 
The paper is organized as follows. In Section~\ref{s:Proof}, we prove Theorems~\ref{p:existence}, 
\ref{p:homogenization}, and 
\ref{prop:OptCondMultEig}.
In Section~\ref{s:disk}, we further analyze the case when $\Omega$ is a disk, proving Theorem~\ref{p:Disk}. 
In Section~\ref{s:CompMeth}, we describe computational methods for weighted Steklov eigenvalues \eqref{e:relax} as well as for approximating the solutions to the non-convex, non-smooth extremal eigenvalue problems  \eqref{e:RelaxMin} and \eqref{e:RelaxMax}. 
In Section~\ref{s:CompRes}, these computational methods are used in  computational experiments for the disk. 
In Appendix~\ref{a:counterexample}, we give a simple example of a Sturm-Liouville operator that has non-convex higher eigenvalues. 

\subsection*{Acknowledgments}
The authors would like to thank the American Institute of Mathematics (AIM) for hosting a SQuaREs workshop on ``Theoretical, Asymptotic, and Numerical Analysis of Extremal Steklov Eigenvalue Problems'', where this project was initiated. The authors would like to thank Weaam A. Alhejali for initial discussions on this project. The authors gratefully acknowledge Seyyed Abbas Mohammadi for valuable discussions that led to the early discovery of numerical evidences of Appendix~\ref{a:counterexample}.


\section{Proofs of Theorems~\ref{p:existence}, \ref{p:homogenization}, and \ref{prop:OptCondMultEig}} \label{s:Proof}
To prove Theorem~\ref{p:existence}, 
we first recall a few basic properties of the admissible set of densities, $\calA_\alpha$, and establish some properties of the eigenvalues $\lambda_k(\rho)$ of the weighted Steklov problem \eqref{e:relax}. 

\begin{lem}[\cite{Friedland_1977}] \label{l:PropA}
Fix $\alpha \in (0, 1]$. The following are properties of the admissible set of densities, $\calA_\alpha$, defined in \eqref{e:Ad}. 
\begin{enumerate}
\item[(a)] $\calA_\alpha$ is a convex subset of $L^\infty(\Gamma)$. 
\item[(b)] Extremal points of $\calA_\alpha$ are of the ``bang-bang'' form: 
$ \rho(s) = \one_S(s)$, $S \subset \Gamma$, $\abs{S} = \alpha\abs{\Gamma}$.  
\item[(c)] $\calA_\alpha$ is compact for the weak* convergence. 
\end{enumerate}
\end{lem}

Let $H^q(\Gamma)$ be the Sobolev space of order $q$, with $H^0 = L^2$, and consider its mean-zero subspace, $H_*^q(\Gamma)\coloneqq \left\{\varphi\in H^q(\Gamma)\colon \int_\Gamma\varphi \, ds = 0\right\}$. Denote by $\Lambda \in \mathcal{B}(H_*^{1/2},H_*^{-1/2})$ the (bounded) Dirichlet-to-Neumann operator for the constant boundary weight $\rho\equiv \one$ and by $M_\rho\in \mathcal{B}(L^2)$ the multiplication operator by boundary weight $\rho\in L^\infty(\Gamma)$. 

\begin{lem} \label{l:PropSigma}
The following are spectral properties of the weighted Steklov eigenproblem \eqref{e:relax}.
\begin{enumerate}
\item[(a)] The operator $T\in \mathcal B (H_*^{-1/2}(\Gamma), H_*^{1/2}(\Gamma))$, defined by $T(\rho)\coloneqq \Lambda^{-1/2}\,M_\rho\,\Lambda^{-1/2}$, is a self-adjoint Hilbert-–Schmidt operator with spectrum given by the nonzero reciprocal weighted Steklov spectrum, $\{\lambda^{-1}_k(\rho)\}_{k\geq 1}$. 
Moreover, let $\{U_j\}_{j\ge0}$ be the harmonic Steklov eigenfunctions in $\Omega$ for the constant density $\rho\equiv \one$, chosen such that $\phi_j = U_j\big|_{\Gamma}$ are orthonormal in $L^2(\Gamma)$, and satisfy $\Lambda \phi_i = \sigma_i \phi_i$, $i \geq 0$ with eigenvalues $0 = \sigma_0 < \sigma_1\le \sigma_2\le \cdots\nearrow\infty$. Define the infinite matrices $D\coloneqq\operatorname{diag}(\sigma_i)_{i\ge 1}$ and $B(\rho)$ with entries given by $B_{ij}(\rho)= \int_\Gamma \rho\phi_i\phi_j\,ds$, $i,j \geq 1$. Then $(\lambda, u)$ is a non-trivial eigenpair of \eqref{e:relax} with 
\begin{equation} \label{e:uBasis}
u(x) = \sum_{k\in \N} c_k \sigma_k^{-\frac{1}{2}} U_k(x) 
\end{equation}
if and only if $(1/\lambda,c)$, with $c = \{c_k\}_{k\in\N}\in\ell^2(\N)$, is an eigenpair of the self-adjoint Hilbert--Schmidt operator $A\in \mathcal B (\ell^2(\N) )$ with matrix representation $A(\rho)\coloneqq D^{-1/2}\,B(\rho)\,D^{-1/2}$.

\item[(b)] The map $\rho\mapsto\lambda_k(\rho)$ is continuous on $\calA_\alpha$ with respect to weak\(^*\) topology on $L^\infty(\Gamma)$. 

\item[(c)] Using notation from (a), the reciprocal principal eigenvalue $ \lambda_1^{-1}(\rho)$ has the variational characterization 
\begin{equation}
\label{e:Var3}
\lambda_1^{-1}(\rho)= \max_{c \in \ell^2(\N), \, (c,b) = 0} \frac{c^t A(\rho) c}{c^t c}, 
\end{equation} 
where $b = \{b_k\}_{k\in\N}\in\ell^2(\N)$ is given by $b_k = \int_{\Gamma} \sigma_k^{-1/2} \phi_k(s) \rho(s) \, ds$. 

\item[(d)] The Fr\'echet derivative of a simple eigenvalue $\lambda\colon L^2(\Gamma) \to \R$ with corresponding normalized eigenfunction, $u$, is given by
$\frac{\delta \lambda}{\delta \rho} = - \lambda u^2$.

\item[(e)] (Gateaux derivative of $\lambda_k(\rho)$) For $\rho$ fixed, suppose that $\lambda_{k_0}$ is a multiplicity $m_{k_0}$ eigenvalue with an (arbitrarily chosen) basis $\{u_{j,k_0}\}_{j \in [m_{k_0}]}$ of the corresponding eigenspace, orthonormal in $L^2(\Gamma; \rho)$. 
For a perturbed density $\rho + \vareps\delta\rho$ with $\delta\rho\in L^\infty(\Gamma)$, the perturbed eigenvalues are given by $\lambda_{k_0} + \vareps \nu_i +\mathcal{O}(\vareps^2)$, $i \in [m_{k_0}]$, where $\nu_i$ are eigenvalues of the matrix $M = M(\delta \rho)\in\R^{m_{k_0} \times m_{k_0}}$ with entries
$ M_{ij} = - \lambda_{k_0} \int_\Gamma u_{i, k_0}(s)u_{j, k_0}(s)\delta\rho(s)\, ds$.   
\end{enumerate}
\end{lem}

\begin{proof}
{\bf (a)} The Rayleigh quotient for a harmonic $u$ with boundary trace $u|_\Gamma=\sum_{j\ge1}y_j\phi_j$ can be written as
$\lambda^{-1} = \frac{\int_\Gamma \rho u^2\,ds}{\int_\Omega \abs{\nabla u}^2 dx}
=\frac{y^* B(\rho)\,y}{y^* D\,y}$,
since $\int_\Omega|\nabla U_i|^2\,dx=\sigma_i$ for the harmonic extension $U_i$ of $\phi_i$.
Thus, the weighted Steklov generalized eigenproblem,
$B(\rho)\,y=\lambda^{-1} \, D\,y$, 
is equivalent to the self-adjoint eigenproblem, $A(\rho)c=\lambda^{-1} \, c$, obtained by the similarity transformation, $c\coloneqq D^{1/2}y$. 
Equivalently, $T(\rho)\in\mathcal B(H_*^{-1/2}(\Gamma),H_*^{1/2}(\Gamma))$ defined by $T(\rho)\coloneqq\Lambda^{-1/2}\,M_\rho\,\Lambda^{-1/2}$, when expressed in the Steklov basis, has matrix entries
\begin{equation}
\label{e:TMatrix}
T_{ij}(\rho) = \frac{1}{\sqrt{\sigma_i\sigma_j}}\int_\Gamma \rho\phi_i\phi_j\,ds = A_{ij}(\rho).
\end{equation}
Hence, the nonzero reciprocals of the weighted Steklov eigenvalues coincide with the nonzero eigenvalues of the compact self-adjoint operator $T(\rho)$.

We now show that $A\colon\ell^2(\N)\to\ell^2(\N)$ is a Hilbert-Schmidt operator. Using $\abs{\int_\Gamma \rho\phi_i\phi_j \, ds}\le \|\rho\|_\infty$ and the Cauchy-Schwarz inequality, we get the bound 
\[
\|A(\rho)\|_\mathrm{HS}^2 = \sum_{i,j\ge 1} \abs{A_{ij}}^2 = 
\sum_{i,j\ge 1} \frac{\abs{\int_\Gamma \rho\phi_i\phi_j \, ds}^2}{\sigma_i \sigma_j} 
\le \|\rho\|_\infty^2\left(\sum_{i\ge 1}\frac{1}{\sigma_i^2}\right)\left(\sum_{j\ge 1} \frac{1}{\sigma_j^2}\right) = \|\rho\|_\infty^2\left(\sum_{k\ge1} \sigma_k^{-2}\right)^2. 
\] 
By Weyl's law for Steklov eigenvalues on smooth planar domains, we have \(\sigma_k\sim C k\) for large \(k\). Consequently, the infinite sum is finite and this shows that \(A(\rho)\) is Hilbert--Schmidt.


{\bf (b)} We follow the proof of \cite[Prop. 2.2]{Cox_2003}. Suppose $\rho_n\ws\rho$ in $L^\infty(\Gamma)$ and  $\lambda_k^{(n)}$ is the $k$-th eigenvalue satisfying \eqref{e:relax} for density $\rho_n$ with an associated eigenfunction $u_k^{(n)}$, normalized so that 
\begin{equation} \label{e:PropEq} 
\int_\Gamma \rho_n\left(u_k^{(n)}\right)^2 ds = 1 \qquad \textrm{ and } \qquad \lambda_k^{(n)} 
= \int_\Omega \abs{\nabla u_k^{(n)}}^2 dx. 
\end{equation}
From Poincar\'e inequality, \eqref{e:PropEq}, and \eqref{e:UpperBnd}, we have that $\left\{u_k^{(n)} \right\}_n$ is bounded in $H^1(\Omega)$: 
\[ \int_\Omega \left(\left(u_k^{(n)}\right)^2 + \abs{\nabla u_k^{(n)}}^2\right) dx \leq (1 + C^2) \int_\Omega \abs{\nabla u_k^{(n)}}^2 dx = (1 + C^2)\lambda_k^{(n)} \leq  (1 + C^2) \frac{8\pi k}{ \alpha\abs{\Gamma}}. \] 
Hence, there exists $u_k \in H^1(\Omega)$ with $u_k^{(n)} \rightharpoonup u_k$ in $H^1(\Omega)$, $u_k^{(n)}  \to u_k$ in $L^2(\Omega)$, and the traces satisfying 
$u_k^{(n)}\big|_\Gamma\to u_k\big|_\Gamma$ in $L^2(\Gamma)$. By \eqref{e:PropEq}, $\lambda_k^{(n)} \to \lambda_k$. Passing to the limit as $n \to \infty$ of the weak form, 
$\int_\Omega \nabla u_k^{(n)}\cdot\nabla v\, dx = \lambda_k^{(n)} \int_\Gamma \rho_n u_k^{(n)}v\, ds$ for all $v\in H^1(\Omega)$, 
we have that 
\[ \int_\Omega \nabla u_k\cdot\nabla v\, dx = \lambda_k \int_\Gamma \rho u_k v\, ds, \qquad \forall v \in H^1(\Omega). \] 
This shows that $\lambda_k$ is the $k$-th eigenvalue for $\rho$ with an associated eigenfunction $u_k$.

{\bf (c)} The variational characterization \eqref{e:Var3} follows from part (a). 

{\bf (d)} Denote the variation with respect to $\rho$ with a dot. Taking variations of $\lambda = \int_\Omega \abs{\nabla u}^2 dx$ and using Green's identity (integration by parts), we obtain
\begin{align*}
\dot\lambda & = 2\int_{\Omega} \nabla u\cdot \nabla\dot u\, dx = 2\int_{\Gamma} (\del_nu)\dot u\, ds - 2\int_\Omega (\Delta u)\dot u\, dx \stackrel{\eqref{e:relax}}{=} 2\lambda\int_{\Gamma} \rho u\dot u\, ds. \end{align*} 
The variation of the normalization condition, $\int_{\Gamma} \rho u^2\, ds = 1$, gives
$ 2\int_\Gamma \rho u \dot u\, ds + \int_{\Gamma} \dot\rho u^2\, ds = 0$. Combining the two equations gives the desired result. 

{\bf (e)} Let $\{\lambda_k\}_{k\in\N}$ be the eigenvalues of \eqref{e:relax} for some fixed density $\rho$, with multiplicities $m_k\in\N$ and corresponding $L^2(\Gamma; \rho)$-orthonormal bases $\{u_{j,k}\}_{j\in [m_k]}$ for their eigenspaces. Let $(\lambda_\vareps, u_\vareps)$ be an eigenpair of \eqref{e:relax} with perturbed density $\rho + \vareps\delta\rho$, for some fixed $\delta\rho\in L^{\infty}(\Gamma)$. Since the eigenbasis for the unperturbed problem is an orthogonal basis for harmonic functions in $\Omega$, for a fixed $k_0\in\N$, we may write a solution $(\lambda_\vareps, u_\vareps)$ as follows: 
\begin{subequations} \label{e:perturbed_ansatz} 
\begin{align}
\label{e:perturbedansatzeigenvalue}
     \lambda_\vareps & = \lambda + \vareps \nu + \mathcal{O}(\vareps^2), \\ 
\label{e:perturbedansatzfunction}
     u_\vareps & = \sum_{k=1}^\infty \sum_{j=1}^{m_k} \Big(\delta_{k,k_0}a_j + \vareps b_{j,k} +\mathcal{O}(\vareps^2)\Big) u_{j,k}, 
\end{align}
\end{subequations} 
where $\delta_{k, k_0}$ denotes the usual Kronecker delta. Inserting \eqref{e:perturbed_ansatz} into \eqref{e:relaxb} with the perturbed density $\rho + \vareps\delta\rho$ and using the boundary condition $\partial_n u_{j,k} = \rho\lambda_k u_{j,k}$ on $\Gamma$, we obtain the following equation on $\Gamma$:
{\small 
\begin{align*}
& \sum_{k=1}^\infty \sum_{j=1}^{m_k} \Big(\delta_{k,k_0}a_j + \vareps b_{j,k} +\mathcal{O}(\vareps^2)\Big)\rho\lambda_k u_{j,k} 
=  \sum_{k=1}^\infty \sum_{j=1}^{m_k} \Big(\rho + \vareps\delta\rho\Big)\Big(\lambda + \vareps\nu +\mathcal{O}(\vareps^2)\Big) \Big(\delta_{k,k_0}a_j + \vareps b_{j,k} +\mathcal{O}(\vareps^2)\Big) u_{j,k}. 
\end{align*}}%
At $\mathcal{O}(1)$, we obtain $\lambda = \lambda_{k_0}$. At $\mathcal{O}(\vareps)$, we have that
\begin{align}
\label{e:orderepsilonperturbed}
    \sum_{k=1}^\infty \sum_{j=1}^{m_k} b_{j,k} \rho (\lambda_k - \lambda_{k_0}) u_{j,k} & = \sum_{j = 1}^{m_{k_0}} \left(\lambda_{k_0}\delta\rho + \nu\rho\right)a_ju_{j, k_0}. 
\end{align}
We now multiply \eqref{e:orderepsilonperturbed} by $u_{i,k_0}$ for $i\in [m_{k_0}]$ fixed and integrate over $\Gamma$. The left hand side is zero by the $\rho$-orthogonality of $\{u_{j,k_0}\}_{j \in [m_{k_0}]}$ on $\Gamma$, and we have
\begin{align*} 
    -\sum_{j=1}^{m_{k_0}} \lambda_{k_0}a_j \int_\Gamma u_{j,k_0}u_{i,k_0}\delta\rho\, ds = \sum_{j=1}^{m_{k_0}} \nu a_j \int_\Gamma u_{j,k_0}u_{i,k_0}\rho\, ds = \nu a_i, \qquad i\in [m_{k_0}], 
\end{align*}
where the last equality is also a consequence of the $\rho$-orthogonality of $\{u_{j,k_0}\}_{j \in [m_{k_0}]}$ on $\Gamma$. In matrix form, this is equivalent to $Ma = \nu a$; the eigenvalues of $M$ give the perturbed eigenvalues and the eigenvectors determine the basis of eigenfunctions to be perturbed. 
\end{proof}

\begin{proof} [Proof of Theorem~\ref{p:existence}.]
Using the direct method of the calculus of variations, the existence of a minimizer in  \eqref{e:RelaxMin} and a maximizer in \eqref{e:RelaxMax} follows from \eqref{e:UpperBnd},  Lemma~\ref{l:PropA}(c), and Lemma~\ref{l:PropSigma}(b). 
\end{proof}

Theorem~\ref{p:homogenization} follows directly from 
Lemma~\ref{l:PropSigma}(b) applied to $\rho_n = \one_{S_n}$.

\begin{proof}[Proof of Theorem~\ref{prop:OptCondMultEig}.]  
Let $\rho_k^\bigtriangledown \in \mathcal{A}_\alpha$ be a local minimizer for \eqref{e:RelaxMin}, let $E_k\subset L^2(\Gamma; \rho_k^{\bigtriangledown})$ denote the $m_k$--dimensional eigenspace corresponding to $\lambda_k = \lambda_k(\rho_k^\bigtriangledown)$, and 
let $\{u_{i, k}\}_{i \in [m_k]} \subset E_k $ be an (arbitrarily chosen) orthonormal basis of $E_k$, normalized so that 
$\int_\Gamma  u_{i, k}(s) u_{j, k}(s) \rho_k^{\bigtriangledown}(s) \, ds = \delta_{ij}$.
By Lemma~\ref{l:PropSigma}(e), for any admissible perturbation $\delta\rho$, i.e.  $\int_\Gamma \delta\rho\, ds =0$ and
$0\le \rho_k^\bigtriangledown +\vareps\delta\rho\le1$ for $|\vareps|\ll 1$, the first-order eigenvalue perturbations are eigenvalues of the symmetric matrix 
\[ M(\delta\rho)\ =\ -\lambda_k \left(\int_\Gamma u_{i, k}(s)u_{j, k}(s)\delta\rho(s)\,ds \right)_{i,j \in [m_k]}. \] 
Since $\rho_k^{\bigtriangledown}$ is a local minimizer, $M(\delta \rho)$ must be positive definite for every admissible $\delta \rho$, i.e., 
\begin{equation} \label{e:LocOpt}
\int_\Gamma q_a(s)\delta\rho(s) \, ds \leq 0,
\qquad \quad 
\forall a \in \R^{m_k}, \, |a| = 1, 
\end{equation}
where $q_a(s)\coloneqq \left(\sum_{i \in [m_k]} a_i u_{i, k}(s)\right)^2$. 
Write $F(s) \coloneqq \sum_{i \in [m_k]} u_{i, k}^2(s)$, $s \in \Gamma$ and 
define the following subsets of $\Gamma$, 
\[ A_0\coloneqq \left\{ \rho_k^{\bigtriangledown} = 0\right\}, \quad 
A_1\coloneqq \left\{ \rho_k^{\bigtriangledown} = 1 \right\}, \quad \textrm{and} \quad 
A\coloneqq \left\{ 0 < \rho_k^{\bigtriangledown} < 1 \right\}. \] 

We first consider the set $A$. For $d \in \left(0,\frac{1}{2}\right)$, define
$A_d\coloneqq \left\{s \in \Gamma \colon  d \leq \rho_k^{\bigtriangledown}(s) \leq 1 -d \right\} \subset A$. 
Fix $d >0$, let $S \subset A_d$ be an arbitrary measurable set, and consider the perturbation 
$\delta\rho_S = \one_S - \frac{\abs{S}}{\abs{A_d}}\one_{A_d}$.
Note that $\delta\rho_S$ and $- \delta \rho_S$ are admissible for sufficiently small $\vareps > 0$. 
Applying \eqref{e:LocOpt} to both $\delta \rho_S$ and $-\delta \rho_S$ gives that 
$\int_{\Gamma} u_{i, k}(s)u_{j, k}(s)\delta\rho_S\, ds = 0$, $i,j \in [m_k]$, which is equivalent to 
\[ \frac{1}{\abs{A_d}} \int_{A_d} u_{i, k}(s)u_{j, k}(s)\, ds  \, = \,  
\frac{1}{\abs{S}} \int_{S} u_{i, k}(s)u_{j, k}(s)\, ds, \qquad i,j \in [m_k]. \] 
Since $S\subset A_d$ is arbitrary, every pairwise product $u_{i, k}(s)u_{j, k}(s)$ is constant on $A_d$. In particular, 
\[ F(s) = \sum_{i \in [m_k]} u_{i, k}^2(s) \equiv c, \qquad s \in A_d, \] 
for some constant $c > 0$.  Using the continuity of $F$ and writing $A = \cup_{d \in \left(0,\frac{1}{2}\right)} A_d$, we get
$$
F(A) = F \left( \cup_{d \in \left(0,\frac{1}{2}\right)} A_d \right) = \cup_{d \in \left(0,\frac{1}{2}\right)} F(A_d) = c.
$$

Next, assume that $\abs{A} >0$ and consider the set $A_0$. Let $S \subseteq A_0$ be a measurable set and consider the perturbation $\delta \rho = \one_S - \frac{\abs{S}}{\abs{A_d}} \one_{A_d}$, for a small fixed $d\in \left(0,\frac{1}{2}\right)$. This perturbation is admissible for sufficiently small $\vareps >0$. 
Applying \eqref{e:LocOpt}, we obtain 
$
\frac{1}{\abs{S}} \int_S q_a(s) \,ds
\leq 
\frac{1}{\abs{A_d}} \int_{A_d} q_a(s) \,ds$.
Choosing $a = \delta_j$ and summing over $j \in [m_k]$, we obtain 
$\frac{1}{\abs{S}} \int_S F(s) \,ds
\leq 
\frac{1}{\abs{A_d}} \int_{A_d} F(s) \,ds = c$. 
Since $S\subseteq A_0$ is arbitrary, we conclude that $F(s) \leq c$ on $A_0$. 
Arguing similarly, assuming $\abs{A} > 0$, we have that $F(s) \geq c$ on $A_1$. This establishes \eqref{e:minRhoCondMult} under the assumption that $\abs{A} > 0$. 

Finally, assume that $\abs{A} = 0$ and consider the sets $A_0$ and $A_1$. 
Fix arbitrary measurable sets $S\subseteq A_0$ and $T \subseteq A_1$ and consider the perturbation 
$\delta \rho = \frac{1}{\abs{S}} \one_S - \frac{1}{\abs{T}} \one_T$, which is admissible for sufficiently small $\vareps > 0$.  From \eqref{e:LocOpt}, we obtain 
$$
\frac{1}{\abs{S}} \int_S \Big(\sum_{i \in [m_k]} a_iu_{i, k}(s)\Big)^2\, ds
\leq 
\frac{1}{\abs{T}} \int_T \Big(\sum_{i \in [m_k]} a_iu_{i, k}(s)\Big)^2\, ds, 
\qquad \quad 
\forall a \in \R^{m_k}, |a| = 1. 
$$
Again, choosing $a = \delta_j$ and summing over $j \in [m_k]$, we obtain 
$\frac{1}{\abs{S}} \int_S F(s) \,ds
\leq 
\frac{1}{\abs{T}} \int_T F(s) \,ds$. 
Taking $c\in \Big[\underset{A_0}{\mathrm{ess\, sup}} \, F, \, \underset{A_1}{\mathrm{ess\, inf}} \, F \Big]$,  we have that $F \leq c$ a.e. on $A_0$ and $F \geq c$ a.e. on $A_1$. Since $\abs{A} = 0$, this establishes \eqref{e:minRhoCondMult}. 

A similar argument shows that any local maximizer $\rho_k^\triangle\in\calA_\alpha$ for  \eqref{e:RelaxMax} satisfies \eqref{e:maxRhoCondMult}. 
\end{proof}

\begin{rem}
One can also prove Theorem~\ref{prop:OptCondMultEig} by formulating the minimization problem \eqref{e:RelaxMin} as the unconstrained problem $\min_\rho \lambda_k(\rho) + \delta_{\calA_\alpha}(\rho)$, where $\delta_{\calA_\alpha}$ is the indicator function. Fermat's theorem gives that a necessary condition for optimality is that 
$0 \in \del f(\rho^\bigtriangledown) + N_{\calA_\alpha}(\rho^\bigtriangledown)$,
where the \emph{normal cone} is given by 
$N_{\calA_\alpha}(\rho) = \{\mu \one + \nu_+ - \nu_- \colon 
\mu \in \R, \, 
\nu_{\pm} \in L_+^\infty, \, 
\textrm{supp}(\nu_+) \subseteq A_1, \, 
\textrm{supp}(\nu_-) \subseteq A_0 \}$.  
Arguing similar to the above proof then gives \eqref{e:minRhoCondMult}.
\end{rem}


\section{Analysis of extremal eigenvalues for a disk} \label{s:disk}
In this section, we consider the case when $\Omega = \D$, so that $\Gamma = \S^1$ and prove Theorem~\ref{p:Disk}. We begin by noting that the constant density $\rho(\theta) = \alpha > 0$ is admissible, since $\int_0^{2\pi} \alpha\, d\theta = 2\pi\alpha = \alpha\abs{\Gamma}$.

\subsection{Local analysis for constant density and proof of Theorem~\ref{p:Disk}(2)} 
In the case of a unit disk $\D$, the weighted Steklov eigenvalues for $\rho = \alpha > 0$ are given by $\lambda_0 = 0$ and $\lambda_{2j - 1} = \lambda_{2j} = j/\alpha$ for $j\in\N$. 
Consider the perturbed density $\rho(\theta) = \alpha + \vareps\delta\rho(\theta)\in [0, 1]$, where $\delta\rho$ satisfies $\int_0^{2\pi} \delta\rho(\theta)\, d\theta = 0$ so that $\alpha + \vareps\delta\rho \in \calA_\alpha$. For fixed $j\geq 1$, by Lemma~\ref{l:PropSigma}(e), the 
first variation of the eigenvalues $\lambda_{2j - 1}$ and $\lambda_{2j}$ are the eigenvalues of the following $2\times 2$ matrix  
\begin{align*}
M(\delta\rho) & = -\frac{j}{\alpha} \int_0^{2\pi} 
    \begin{pmatrix} 
        u_{2j - 1}^{2}(\theta) & u_{2j - 1}(\theta)u_{2j}(\theta) \\ 
        u_{2j - 1}(\theta)u_{2j}(\theta) &  u_{2j}^{2}(\theta) 
    \end{pmatrix} \delta\rho(\theta)\, d\theta, 
\end{align*} 
where $u_{2j - 1}(\theta) = \sqrt{\frac{1}{\pi\alpha}}\cos(j\theta)$ and $u_{2j}(\theta) = \sqrt{\frac{1}{\pi\alpha}} \sin(j\theta)$. After simplification by using the double angle formula for trigonometric functions, we have 
\[ M(\delta\rho) = -\frac{j}{2\pi\alpha^2} \int_0^{2\pi} 
     \begin{pmatrix} 
         \cos(2j\theta) & \sin(2j\theta) \\ 
         \sin(2j\theta) & -\cos(2j\theta) 
     \end{pmatrix} \delta\rho(\theta)\, d\theta. \] 
Note that the trace of $M$ is zero and the determinant of $M$ is 
\[ \det(M(\delta\rho)) = -\left(\frac{j}{2\pi\alpha^{2}}\right)^{2} \left\{\left(\int_0^{2\pi} \cos(2j\theta)\delta\rho(\theta)\, d\theta\right)^{2} + \left(\int_0^{2\pi} \sin(2j\theta)\delta\rho(\theta)\, d\theta\right)^{2}\right\}. \] 
Expanding $\delta\rho$ as the Fourier series, 
$ \delta\rho(\theta) = \sum_{\ell = 1}^{\infty} \Big(a_\ell\cos(\ell\theta) + b_\ell\sin(\ell\theta)\Big)$, 
we have 
$\det(M) = -\left(\frac{j}{2\alpha^{2}}\right)^{2}\left\{(a_{2j})^{2} + (b_{2j})^{2}\right\}$. Whenever $a_{2j}$ or $b_{2j}$ are not identically zero, we have
$\det(M) < 0$ and the perturbed eigenvalues are given by 
\begin{alignat*}{2}
\lambda_{2j - 1} & = \frac{j}{\alpha} - |\vareps|\frac{j}{2\alpha^{2}}\sqrt{(a_{2j})^{2} + (b_{2j})^{2}} + \mathcal{O}(\vareps^2) && \le \frac{j}{\alpha}, \\ 
\lambda_{2j} & = \frac{j}{\alpha} + |\vareps| \frac{j}{2\alpha^{2}}\sqrt{(a_{2j})^{2} + (b_{2j})^{2}} + \mathcal{O}(\vareps^2) && \ge \frac{j}{\alpha}, 
\end{alignat*}
where $\vareps\to 0$. 
We conclude the following. 
\begin{enumerate} 
\item There is no perturbation $\delta\rho$ that can increase $\lambda_{2j - 1}$ at $\mathcal{O}(\vareps)$. Thus, the constant density $\rho =\alpha$ is a critical point for \eqref{e:RelaxMax} for odd $k\geq 1$. 

\item The admissible perturbation, $\delta\rho(\theta) = a_{2j} \cos(2j\theta)$ for $a_{2j}\neq 0$, will decrease $\lambda_{2j-1}$. Thus, the constant density $\rho = \alpha$ is not a local minimizer for \eqref{e:RelaxMin} for odd $k\geq 1$.

\item There is no perturbation $\delta\rho$ that can decrease $\lambda_{2j}$ at $\mathcal{O}(\vareps)$. Thus, the constant density $\rho =\alpha$ is a critical point for \eqref{e:RelaxMin} for even $k\geq 2$. 
 
\item The admissible perturbation, $\delta\rho(\theta) = a_{2j} \cos(2j\theta)$ for $a_{2j}\neq 0$, will increase $\lambda_{2j}$. Thus, the constant density $\rho = \alpha$ is not a local maximizer for \eqref{e:RelaxMax} for even $k\geq 2$. 
\end{enumerate}

\subsection{Asymptotic results for \texorpdfstring{$\alpha\approx 1$}{alpha approximately 1} and proof of Theorem~\ref{p:Disk}(3)}
Given a fixed $\vareps \ll 1$, we consider the extremal problems \eqref{e:Min} and \eqref{e:Max} for $\alpha = 1- \vareps$, i.e., $\abs{\Gamma_N} = \vareps\abs{\Gamma}$. We first prove the following lemma that bounds the behavior of the perturbed eigenvalues. 

\begin{lem} \label{l:BangBangPert}
Consider the S-N eigenproblem \eqref{e:Steklov} with $\Omega = \D$ and $\Gamma_N = S\subset\S^1$, where the subset $S$ has the form 
$S = \bigcup_{j\in J} S_j = \bigcup_{j\in J} \left(\theta_j - \frac{\abs{S_j}}{2}, \  \theta_j + \frac{\abs{S_j}}{2}\right)$. 
Here, we assume $\sum_{j\in J} \abs{S_j} = 2\pi\vareps>0$ and use the notation $\theta_j$ to denote the midpoint of $S_j$, $j\in J$. For $m = 1, 2, \dots$, the corresponding S-N eigenvalues are given by 
\begin{align*}
    \sigma_{2m - 1} & = m + \nu_{\min}(M) \vareps + \mathcal{O}(\vareps^2), \\ 
    \sigma_{2m} & = m + \nu_{\max}(M)\vareps + \mathcal{O}(\vareps^2), 
\end{align*} 
where $\nu_{\min}$ and $\nu_{\max}$ are the eigenvalues of the $2\times 2$ symmetric matrix $M = M(S)$, given by 
\begin{equation} \label{e:M} 
M = mI + \frac{1}{2 \pi \vareps}\sum_{j\in J} \sin(m\abs{S_j}) 
    \begin{pmatrix}
        \cos(2m\theta_j) & \sin(2m\theta_j) \\ 
        \sin(2m\theta_j) & -\cos(2m\theta_j) 
    \end{pmatrix}. 
\end{equation} 
Moreover, we have that $\nu_{\min}(M)\in [0, m]$ and $\nu_{\max}(M)\in [m, 2m]$, with $\nu_{\min}(M) + \nu_{\max}(M) = 2m$.  
\end{lem} 

\begin{proof}
Since we are perturbing the Steklov spectrum on a disk, we consider the multiplicity 2 eigenvalue $\sigma = m$ with corresponding normalized Steklov eigenfunctions  
$ u_1(r=1,\theta)  = \frac{1}{\sqrt{\pi}} \cos(m\theta)$ and $u_2(r=1,\theta) = \frac{1}{\sqrt{\pi}} \sin(m\theta)$.  
To compute the perturbed eigenvalues, we use Lemma \ref{l:PropSigma}(e) with $\delta \rho(\theta) = - \frac{1}{\vareps}\one_S$. This yields 
\begin{align*}
M & = \frac{m}{\pi \vareps} \int_S 
	\begin{pmatrix} 
		\cos^2(m\theta) & \cos(m\theta)\sin(m\theta) \\ 
		\cos(m\theta)\sin(m\theta) & \sin^2(m\theta) 
	\end{pmatrix} d\theta \\ 
& = \frac{m}{2\pi \vareps}\sum_{j\in J}  
	\begin{pmatrix} 
		\abs{S_j} + \frac{1}{m} \cos(2m\theta_j)\sin(m\abs{S_j}) &  \frac{1}{m}\sin(2m\theta_j)\sin(m\abs{S_j}) \\  
		\frac{1}{m}\sin(2m\theta_j)\sin(m\abs{S_j}) & \abs{S_j} - \frac{1}{m}\cos(2m\theta_j) \sin(m\abs{S_j})
	\end{pmatrix} \\ 
& = mI + \frac{T}{2\pi \vareps} \sum_{j\in J} \frac{\sin(m\abs{S_j})}{T} 
	\begin{pmatrix} 
		\cos(2m\theta_j) & \sin(2m\theta_j) \\  
		\sin(2m\theta_j) & -\cos(2m\theta_j) 
	\end{pmatrix}, 
\end{align*}
where $T\coloneqq \sum_{j\in J} \sin(m\abs{S_j})$ and we use $\sum_{j\in J} \abs{S_j} = 2\pi\vareps$ to obtain $mI$. Recall that the smallest eigenvalue of a symmetric matrix is concave. Together with Jensen's inequality, we find
{\small
\begin{align*}
\nu_{\min}\left(\sum_{j\in J} \frac{\sin(m\abs{S_j})}{T} 
    \begin{pmatrix} 
        \cos(2m\theta_j) & \sin(2m\theta_j) \\  
        \sin(2m\theta_j) & -\cos(2m\theta_j) 
    \end{pmatrix}\right) 
& \ge \sum_{j\in J} \frac{\sin(m\abs{S_j})}{T}\nu_{\min} 
    \begin{pmatrix} 
	\cos(2m\theta_j) & \sin(2m\theta_j) \\ 
	\sin(2m\theta_j) & -\cos(2m\theta_j) 
    \end{pmatrix} \\ 
& = -\sum_{j\in J} \frac{\sin(m\abs{S_j})}{T} 
 \ge -\sum_{j\in J} \frac{m\abs{S_j}}{T} = -\frac{2\pi m\vareps}{T},
\end{align*}}%
where we used the inequality $\sin(z)\le z$ for $0\le z \le \pi$. We conclude that $\nu_{\min}(M) \ge m  - m = 0$. Moreover, with $\tr(M) = 2 m $, we have that $\nu_{\max}(M)\le 2m$. 
\end{proof}

\begin{proof}[Proof of Theorem~\ref{p:Disk}(3).] 
We find configurations of $\Gamma_N$ that saturate the inequalities for $\nu_{\min}(M)$ and $\nu_{\max}(M)$ in Lemma~\ref{l:BangBangPert}. Below, we let $m\in\N$ and split the proof into three cases. 

\vspace{2mm}

\noindent \emph{Case A: \eqref{e:Min} for odd $k\ge 1$ and \eqref{e:Max} for even $k\ge 2$.} 
Consider $\sigma_{2m - 1}$ and $\sigma_{2m}$ for $\Gamma_N = \calI_\vareps^{2m}\subset\S^1$, as defined in Definition~\ref{d:Sn}. Then $\abs{J} = 2m$, $\theta_j = \frac{j\pi}{m}$, and $\abs{S_j} = \frac{\pi\vareps}{m}$ for $j\in [2m]$. For this choice, we have that 
\[ \begin{pmatrix} 
        \cos(2m\theta_j) & \sin(2m\theta_j) \\  
        \sin(2m\theta_j) & -\cos(2m\theta_j) 
    \end{pmatrix} = 
    \begin{pmatrix} 
        1 & 0 \\ 
        0 & -1 
    \end{pmatrix}. \] 
Together with $\sin(m\abs{S_j}) = \sin(\pi\vareps)\approx \pi\vareps$ for $\vareps\ll 1$, we deduce from \eqref{e:M} that $\nu_{\min}(M) = m - m = 0$ and $\nu_{\max}(M) = m + m = 2m$. This implies that for $\vareps\ll 1$, $\Gamma_N = \calI_\vareps^{2m}$ (or $\Gamma_S = \calI_{1 - \vareps}^{2m}$) is a local minimizer for $\sigma_{2m - 1}$ and a local maximizer for $\sigma_{2m}$. 

\vspace{2mm} 

\noindent \emph{Case B: \eqref{e:Max} for odd $k\ge 1$.} Consider $\sigma_{2m - 1}$ and $\sigma_{2m}$ for $\Gamma_N = \calI_\vareps^{2m - 1}\subset\S^1$. Then $\abs{J} = 2m - 1$, $\theta_j = \frac{2j\pi}{2m - 1}$, and $\abs{S_j} = \frac{2\pi\vareps}{2m - 1}$ for $j\in [2m - 1]$. For this choice, we have, from \eqref{e:M}, 
\[ \sum_{j\in [2m - 1]} \sin\left(\frac{2\pi m\vareps}{2m - 1}\right) 
    \begin{pmatrix} 
        \cos(2m\theta_j) & \sin(2m\theta_j) \\  
        \sin(2m\theta_j) & -\cos(2m\theta_j) 
    \end{pmatrix} = \begin{pmatrix} 0 & 0 \\ 0 & 0 \end{pmatrix}, \]
so that $\nu_{\min}(M) = \nu_{\max}(M) = m$. This implies that for $\vareps\ll 1$, $\Gamma_N = \calI_\vareps^{2m - 1}$ (or $\Gamma_S = \calI_{1 - \vareps}^{2m - 1}$) is a local maximizer for $\sigma_{2m - 1}$. 

\vspace{2mm} 

\noindent \emph{Case C: \eqref{e:Min} for even $k\ge 2$.} Consider $\sigma_{2m - 1}$ and $\sigma_{2m}$ for $\Gamma = \calI_\vareps^{2m + 1}\subset\S^1$. Then $\abs{J} = 2m + 1$, $\theta_j = \frac{2j\pi}{2m + 1}$, and $\abs{S_j} = \frac{2\pi\vareps}{2m + 1}$ for $j\in [2m + 1]$. For this choice, we have, from \eqref{e:M}, 
\[ \sum_{j\in [2m + 1]} \sin\left(\frac{2\pi m\vareps}{2m + 1}\right)
    \begin{pmatrix} 
        \cos(2m\theta_j) & \sin(2m\theta_j) \\  
        \sin(2m\theta_j) & -\cos(2m\theta_j) 
    \end{pmatrix} = \begin{pmatrix} 0 & 0 \\ 0 & 0 \end{pmatrix}, \]
so that $\nu_{\min}(M) = \nu_{\max}(M) = m$. This implies that for $\vareps\ll 1$, $\Gamma_N = \calI_\vareps^{2m + 1}$ (or $\Gamma_S = \calI_{1 - \vareps}^{2m + 1}$) is a local minimizer for $\sigma_{2m}$. 
\end{proof}

\section{Computational Methods} \label{s:CompMeth} 
\subsection{Computation of weighted Steklov eigenvalues}
To solve the weighted Steklov eigenproblem \eqref{e:relax} computationally, we first reformulate it using layer potential methods and then discretize the resulting integral equation to obtain a matrix eigenvalue problem. 
In this section, we additionally assume that the domain $\Omega$ is simply connected with $C^2$ boundary. 
The eigenfunction $u(x)$ is represented using a modified single layer potential, $\varphi$, 
\begin{equation} 
\label{steklov_single_layer_modified}
u(x) = \int_{\Gamma} \Phi (x-y)(\varphi(y) - \overline{\varphi})\, ds(y) + \overline{\varphi},
\end{equation}
where $\Phi(x) = \frac{1}{2\pi}\log\abs{x}$ and  $\overline{\varphi}  = \int_{\Gamma}  \varphi(y)\, ds(y) /\int_{\Gamma} \, ds(y)$. 
The modification ensures uniqueness of the solution. 
Taking into account well-known expressions (see {\it e.g.}, \cite{kress1999}) for the
jump of the single layer potential and its normal derivative across
$\Gamma$, the weighted Steklov eigenproblem~\eqref{e:relax} 
reduces to the integral eigenvalue equation for $(\lambda,\varphi)$, 
\begin{equation}
\label{steklov_integral_equation_system}
A[\varphi] = \lambda \rho  B[\varphi]. 
\end{equation}
Here, the boundary operators $A$ and $B$ are defined as 
\begin{align*}
A[\varphi](x) & \coloneqq \int_{\Gamma} \frac{\del\Phi(x - y)}{\del n(x)}(\varphi(y) - \overline{\varphi})\, ds(y) + \frac{1}{2}(\varphi(x) - \overline{\varphi}), \\
B[\varphi](x) & \coloneqq \int_{\Gamma} \Phi (x - y)(\varphi(y) - \overline{\varphi})\, ds(y)+ \overline{\varphi}. 
\end{align*}
Discretizing this equation using the spectrally accurate Nystr\"om quadrature, we obtain the generalized matrix eigenvalue problem, 
\[  {\tt A} {\tt X} =   {\tt R}  {\tt B} {\tt X} {\tt \Lambda}, \] 
where $\tt A$, $\tt B$, and $\tt X$ are discretizations of $A$, $B$, and $\varphi$, respectively. The diagonal matrix $\tt R$ is the discretizations of $\rho$. More details about this method can be found in \cite{Akhmetgaliyev2016}. We implement these methods in MATLAB using \verb+MPSpack+ \cite{mpspack}. An accuracy test in  \cite{Akhmetgaliyev2016} demonstrated that an error less than $10^{-6}$ could be achieved with 512 points for a constant density. 

\subsection{Computation of extremal Steklov eigenvalues}
To approximate solutions to the extremal weighted Steklov eigenvalue problems \eqref{e:RelaxMin} and \eqref{e:RelaxMax}, we apply gradient-based optimization methods. To handle non-simplicity of eigenvalues, we (trivially) reformulate \eqref{e:RelaxMin} as the minimax problem, 
\[ \min_{\rho \in \calA_\alpha}  \ \max_{j\in [k]} \ \lambda_j(\rho). \] 
This minimax problem can be numerically solved  using a sequential quadratic programming (SQP) method \cite{brayton1979new}, as implemented in MATLAB's \verb+fminimax+ function. The gradient of the (numerically simple) eigenvalues is provided using Lemma~\ref{l:PropSigma}(d). We use as convergence criteria a tolerance on either the relative change on stepsize or function value. The maximization problem \eqref{e:RelaxMax} is analogously reformulated as a maximin problem and solved using the same methods.

In this paper, we focus on the numerical implementation of the first few eigenvalues, e.g., $k\le 12$. For moderately large values of $k$, the discretization must be refined to resolve high frequency rapid oscillations. In addition, techniques such as tracking only a subset of eigenvalues and using continuation of subspace methods may be more efficient \cite{men2010bandgap,kangal2018subspace, kao2025semi}. For robustness, we compute all eigenpairs with $j\in [k]$ in our numerical simulations so we can find the corresponding orthonormal eigenfunctions. As there are local optimizers, the algorithm may converge to nearby optimizers which are close to the initial guess. In this paper, we show the best optimizers only.

\section{Computational experiments} \label{s:CompRes}
In this section, we present the results of several computational experiments to illustrate the main results of the paper. In Sections~\ref{s:DiskMin} and \ref{s:DiskMax}, we focus on the minimization problem \eqref{e:RelaxMin} and maximization problem  \eqref{e:RelaxMax} in the case where $\Omega = \D$, respectively. These results provide supporting evidence for Conjecture~\ref{c:Disk}. 

\subsection{Minimizers for the disk} \label{s:DiskMin}
For $\Omega = \D$ and all values of $k$ that we tried, we observe that the minimizer for $\lambda_k$ in \eqref{e:RelaxMin} is the indicator function $\rho_k^\bigtriangledown = \one_{\calI_\alpha^{k+1}}$, where $\calI_\alpha^{k+1}\subset \S^1$ is given in Definition~\ref{d:Sn}. 
Thus, we conjecture
\[ \sigma_k(\calI_\alpha^{k+1}) \ = \ \sigma^\bigtriangledown_{k,\alpha}  \ = \ \lambda^\bigtriangledown_{k,\alpha}, 
\qquad \qquad \qquad k\geq 1, \ \alpha \in (0,1]. \] 
Interestingly, $\rho_k^\bigtriangledown$ is a bang-bang function for all $k\geq 1$. 

\begin{figure}
\centering
$k=1$ \hspace{6.3cm} $k=3$ \hspace{0.8cm}\phantom{x} \\
\includegraphics[trim=20 10 20 10,clip,width=0.45\textwidth]{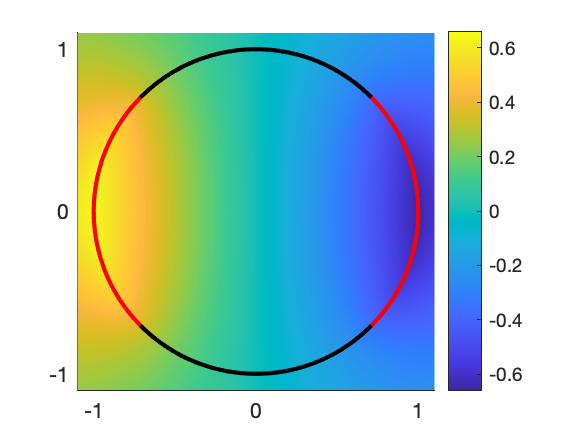}
\includegraphics[trim=20 10 20 10,clip,width=0.45\textwidth]{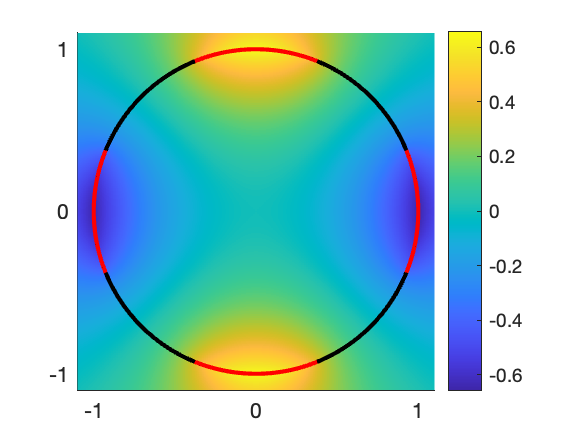}\\
$k=5$ \hspace{6.3cm} $k=7$ \hspace{0.8cm}\phantom{x} \\
\includegraphics[trim=20 10 20 10,clip,width=0.45\textwidth]{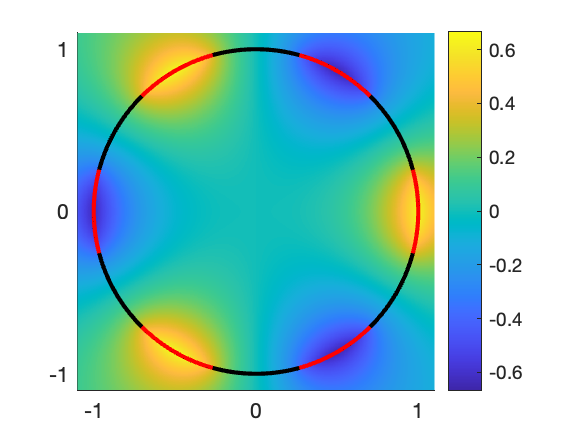}
\includegraphics[trim=20 10 20 10,clip,width=0.45\textwidth]{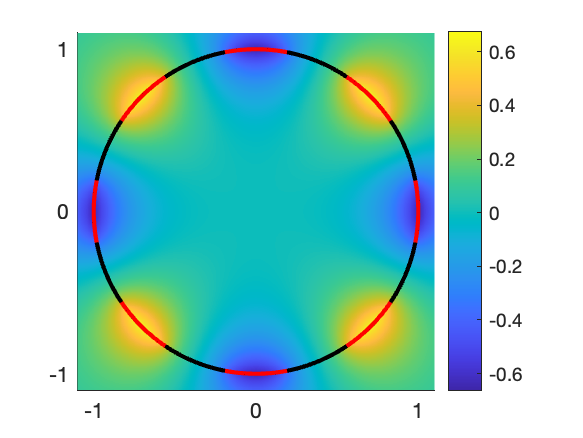}\\
$k=9$ \hspace{6.3cm} $k=11$ \hspace{0.8cm}\phantom{x} \\
\includegraphics[trim=20 10 20 10,clip,width=0.45\textwidth]{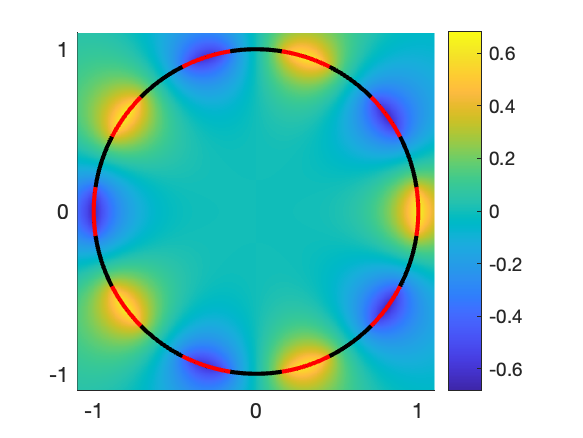}
\includegraphics[trim=20 10 20 10,clip,width=0.45\textwidth]{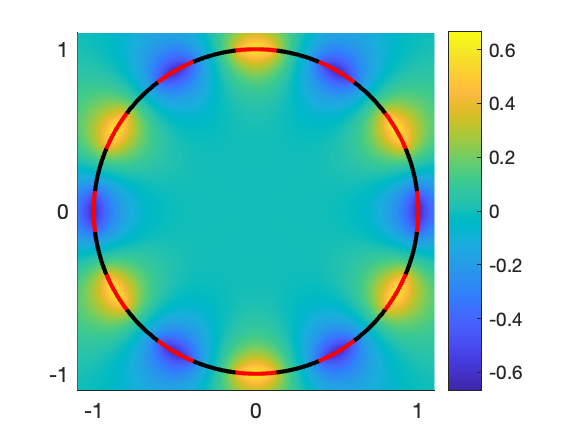}
\caption{Minimizing bang-bang densities, $\rho_k^\bigtriangledown$, and associated eigenfunctions for \eqref{e:RelaxMin} with $\Omega = \D$, $\alpha = 0.5$, and odd $k = 1, 3, 5, 7, 9, 11$ (values indicated). The red part of the boundary is where $\rho=1$, {\it i.e.}, where the Steklov boundary condition is imposed. See Section~\ref{s:DiskMin}.}
\label{f:kOdd} 
\end{figure}

For odd $k$, $\sigma^\bigtriangledown_{k,\alpha}$ have multiplicity one. 
Figure~\ref{f:kOdd} shows the optimal  $\rho_k^\bigtriangledown$ and associated eigenfunction for $\alpha = 0.5$ and $k = 1, 3, 5, 7, 9, 11$. 
The optimal values are given by  
\[ \sigma_{k,\alpha}^\bigtriangledown \ = \ \frac{k+1}{2} \sigma_{1,\alpha}^\bigtriangledown, \qquad \qquad \qquad k\geq 1 \textrm{ odd}, \ \alpha \in (0,1]. \] 

\begin{figure}
\centering
$k=2$\\
\includegraphics[trim=20 10 20 10,clip,width=0.45\textwidth]{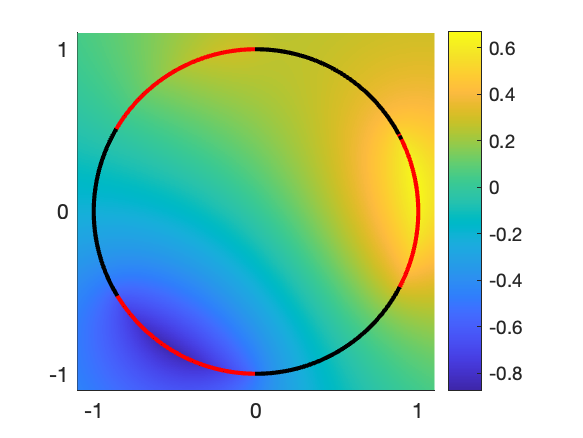}
\includegraphics[trim=20 10 20 10,clip,width=0.45\textwidth]{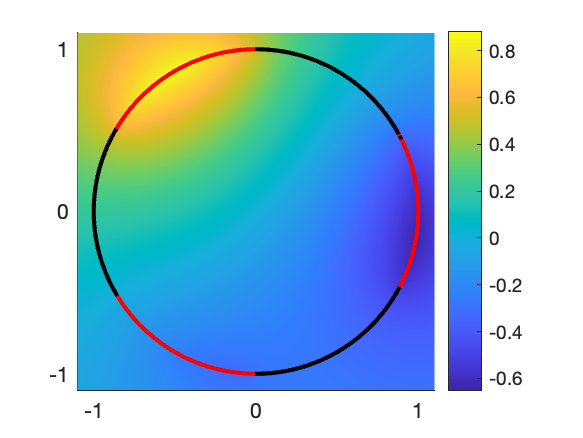}\\
$k=4$\\
\includegraphics[trim=20 10 20 10,clip,width=0.45\textwidth]{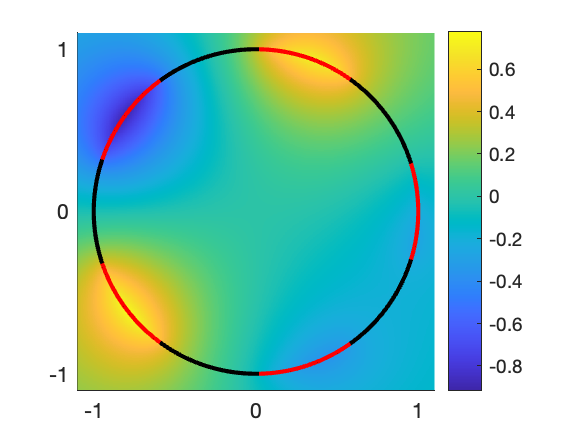}
\includegraphics[trim=20 10 20 10,clip,width=0.45\textwidth]{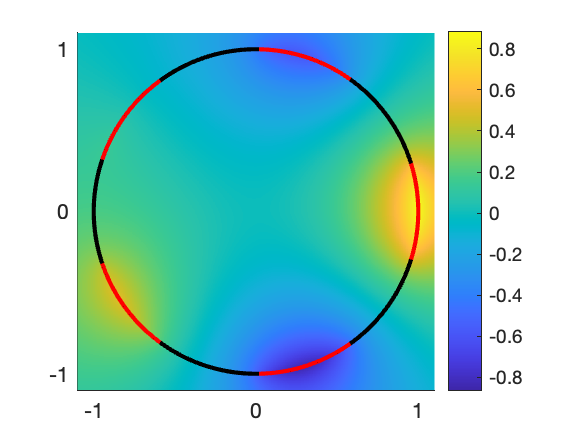}\\
$k=12$\\
\includegraphics[trim=20 10 20 10,clip,width=0.45\textwidth]{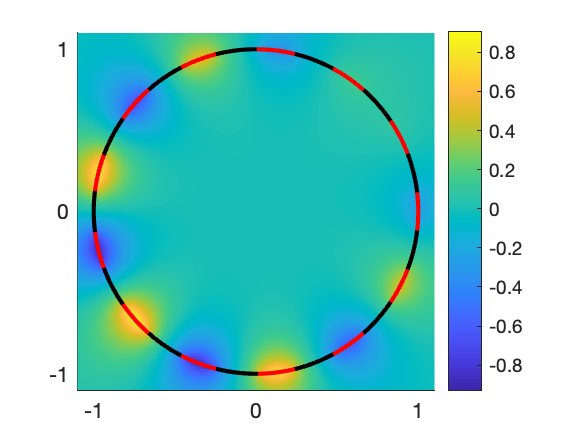}
\includegraphics[trim=20 10 20 10,clip,width=0.45\textwidth]{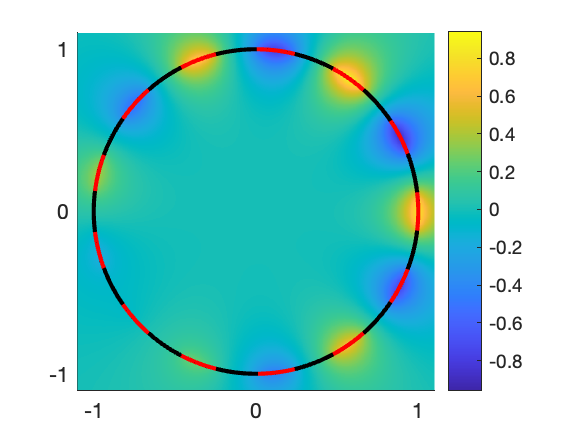}\\
\caption{Minimizing bang-bang densities, $\rho_k^\bigtriangledown$, and associated eigenfunctions for \eqref{e:RelaxMin} with $\Omega = \D$, $\alpha = 0.5$, and even $k=2,4,12$ (values indicated). The red part of the boundary is where $\rho=1$, {\it i.e.}, where the Steklov boundary condition is imposed. Each eigenvalue has multiplicity 2 and a basis for the eigenspace is plotted. See Section~\ref{s:DiskMin}.}
\label{f:kEven} 
\end{figure}

For even $k$, $\sigma^\bigtriangledown_{k,\alpha}$ have multiplicity two. Figure~\ref{f:kEven} shows the optimal $\rho_k^\bigtriangledown$ for $\alpha = 0.5$ and $k = 2, 4,  12$.
The optimal values satisfy 
$ \sigma_{k,\alpha}^\bigtriangledown \to \frac{k+1}{2} \sigma_{1,\alpha}^\bigtriangledown$, $k\geq 2$ even, $\alpha \in (0,1]$.

\begin{figure}
\centering
$k=1$ \\
\includegraphics[width=0.33\textwidth]{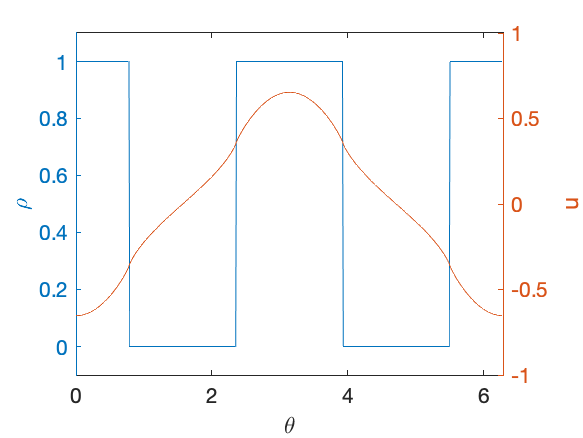}
\includegraphics[width=0.33\textwidth]{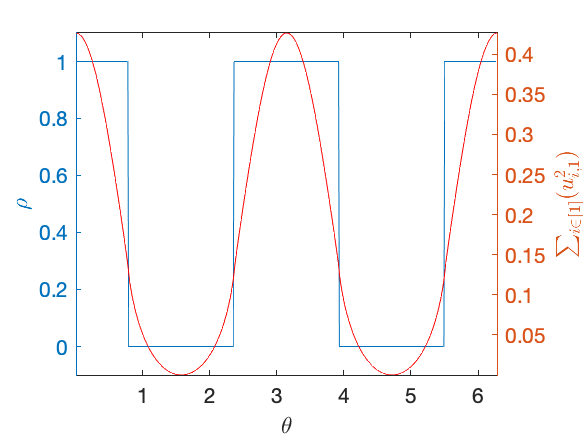}\\
$k=2$ \\
\includegraphics[width=0.33\textwidth]{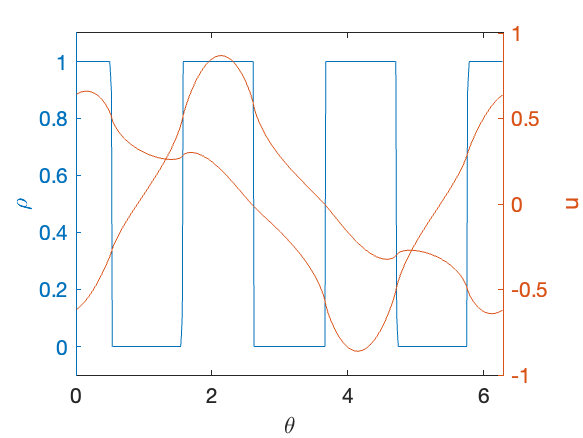}
\includegraphics[width=0.33\textwidth]{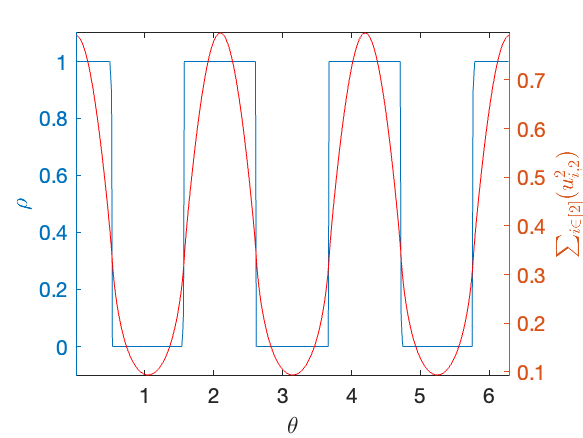}\\
$k=3$ \\
\includegraphics[width=0.33\textwidth]{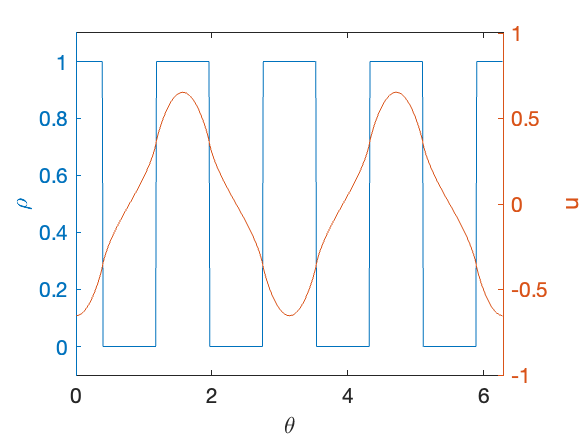}
\includegraphics[width=0.33\textwidth]{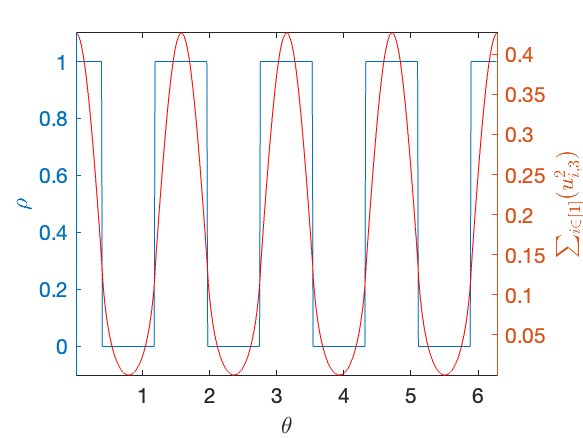}\\
$k=4$ \\
\includegraphics[width=0.33\textwidth]{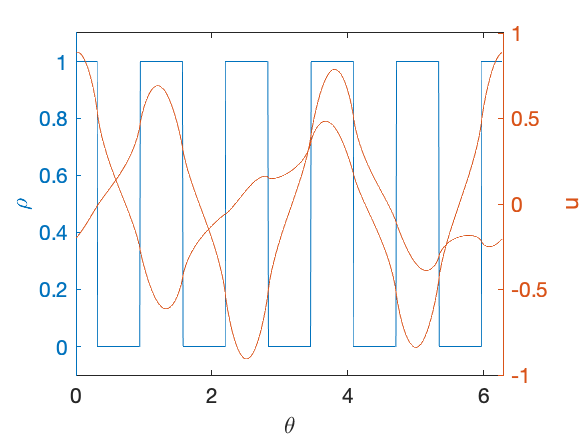}
\includegraphics[width=0.33\textwidth]{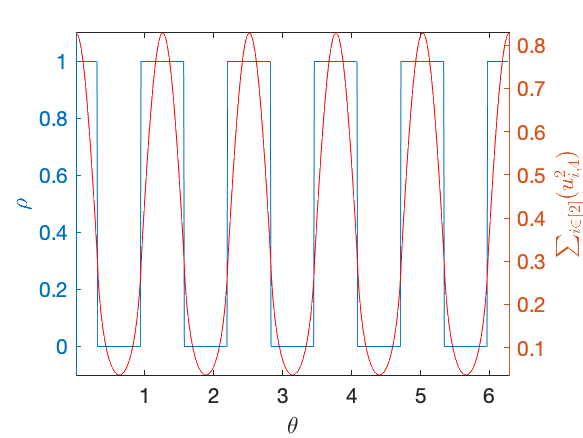}
\caption{Solutions for \eqref{e:RelaxMin} with $\Omega = \D$, $\alpha = 0.5$, and $k = 1,2,3,4$ (values indicated). 
In the left hand panels, the minimizing bang-bang density $\rho_k^\bigtriangledown$ is plotted in blue and the corresponding eigenfunction(s) are plotted in red.
For $k$ odd, the eigenvalue is simple and the right hand panel shows that the optimality condition \eqref{e:minRhoCondMult} is satisfied. 
For $k$ even, the eigenvalue has multiplicity $2$ and the right hand panel shows that the optimality condition \eqref{e:minRhoCondMult} is satisfied for the two eigenfunctions plotted in the left hand panel. See Section~\ref{s:DiskMin}.}
\label{f: min optimiality} 
\end{figure}

In Figure~\ref{f: min optimiality}, we further investigate the optimality conditions established in Theorem~\ref{prop:OptCondMultEig}. 
For odd $k = 1, 3$, $\lambda^\bigtriangledown_{k,\alpha}$ is simple and the eigenfunction satisfies the optimality condition \eqref{e:minRhoCondMult}. Indeed, the optimal $\rho_k^\bigtriangledown$ is bang-bang and the switching occurs at locations for which $u_k^2 \Big|_{\Gamma} (\theta) = c$, for some $c>0$.  
For even $k = 2, 4$, $\lambda^\bigtriangledown_{k,\alpha}$ has multiplicity two and the eigenfunctions satisfy the optimality condition \eqref{e:minRhoCondMult}. Indeed, the optimal $\rho_k^\bigtriangledown$ is bang-bang  and the switching occurs at locations for which $\sum_{i \in [2]} u^2_{i,k} \Big|_{\Gamma} (\theta) = c$, for some $c>0$. 
Eigenfunctions $u_{i,k}$ are plotted in the left panels of Figure~\ref{f: min optimiality}.

\subsection{Maximizers for the disk} \label{s:DiskMax}
For $\Omega = \D$ and $k = 1$, as established in Theorem~\ref{p:Disk}, we numerically observe the optimal density is constant; see Figure~\ref{f:SmallkMaximizer}(top). A basis for the eigenspace (restricted to $\S^1$) is $\{\sin\theta, \cos \theta \}$ and, of course,  $\sin^2\theta + \cos^2 \theta =1$, so that the optimality condition \eqref{e:maxRhoCondMult} is satisfied. Thus, 
$\lambda_{1,\alpha}^\triangle = \frac{1}{\alpha}$. 

\begin{figure}[p]
\centering
$k=2$ \\
\includegraphics[width=0.36\textwidth]{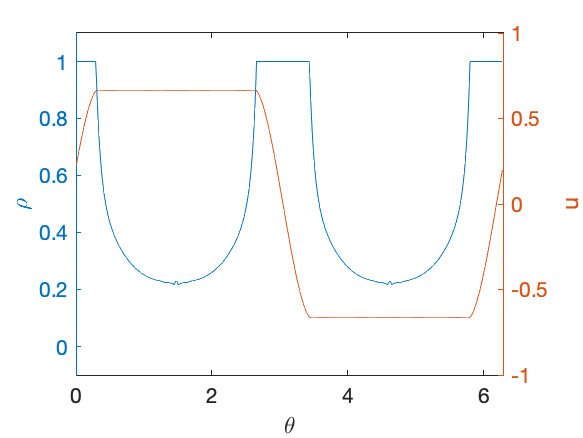}
\includegraphics[width=0.36\textwidth]{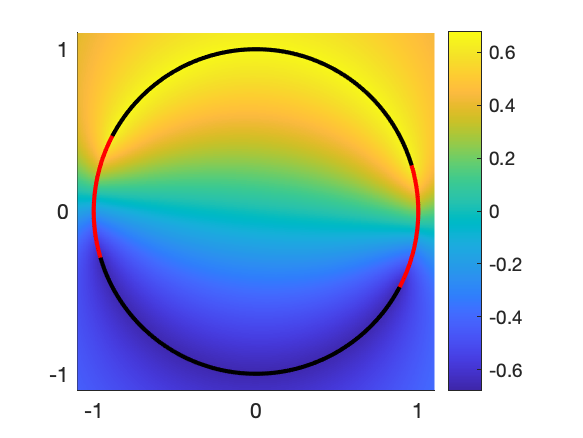}\\
$k=4$ \\
\includegraphics[width=0.36\textwidth]{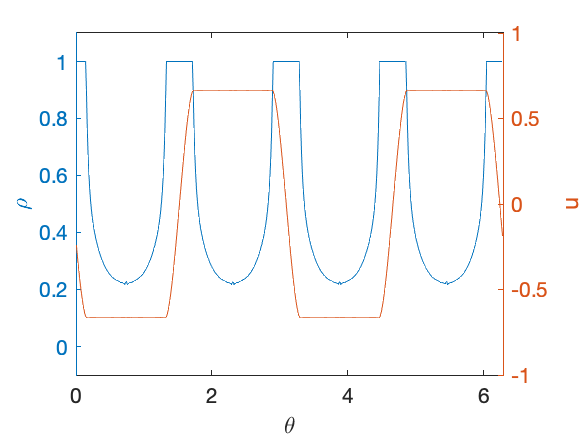}
\includegraphics[width=0.36\textwidth]{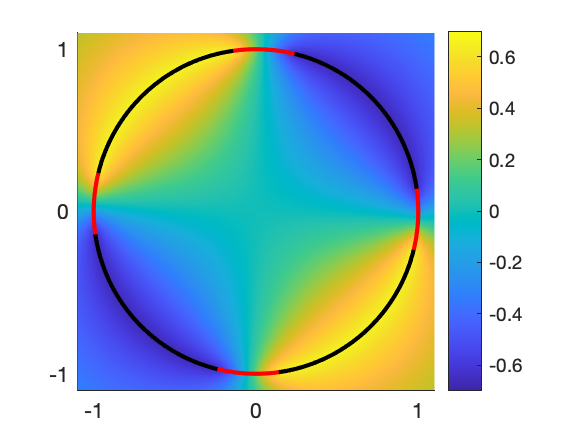}\\
$k=6$ \\
\includegraphics[width=0.36\textwidth]{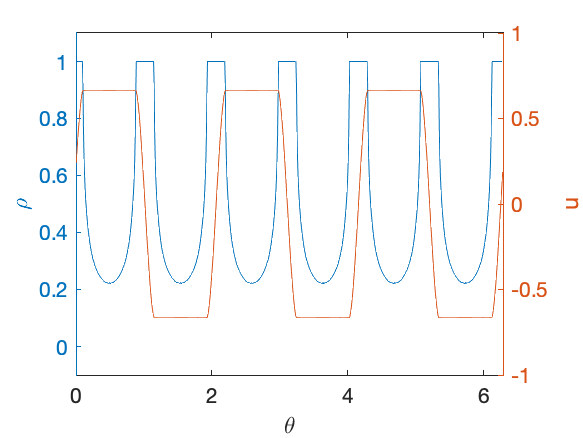}
\includegraphics[width=0.36\textwidth]{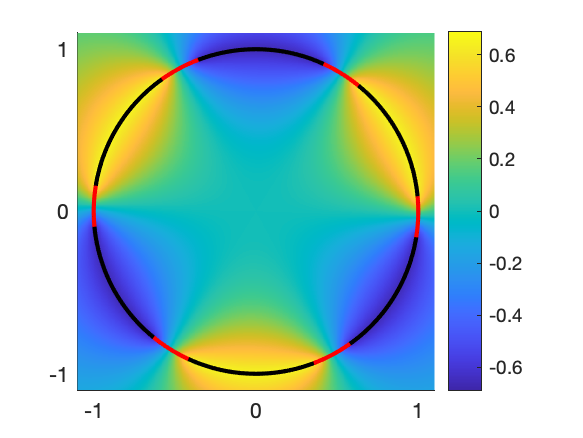}\\
$k=8$ \\
\includegraphics[width=0.36\textwidth]{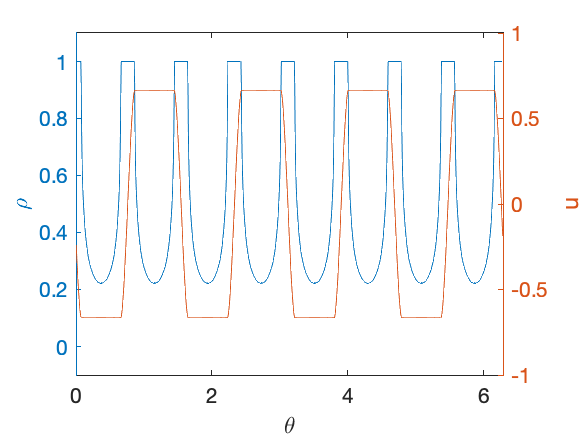}
\includegraphics[width=0.36\textwidth]{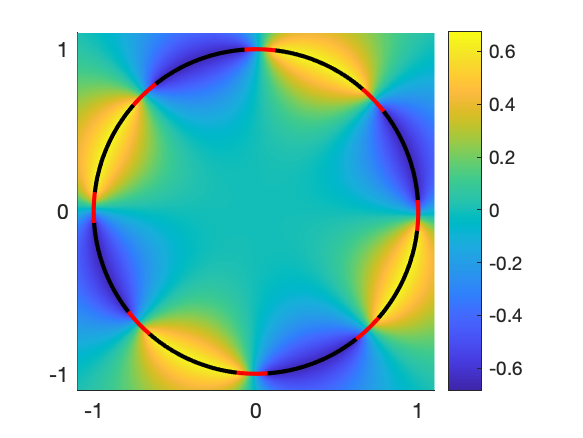}
\caption{{\bf (Left)} Maximizing densities, $\rho_k^\triangle$, and associated eigenfunctions for \eqref{e:RelaxMax} with $\Omega = \D$, $\alpha = 0.5$, and even $k=2,4,6,8$ (values indicated).  
Interestingly, the eigenfunctions are constant on intervals where $\rho \neq 1$. 
{\bf (Right)} The eigenfunctions are plotted on $\Omega =\D$. The red part of the boundary is where $\rho=1$, {\it i.e.}, where the Steklov boundary condition is imposed. See Section~\ref{s:DiskMax}.}
\label{f:EvenkMaximizer}
\end{figure}

\begin{figure}[p]
\centering
$k=1$ \\
\includegraphics[width=0.34\textwidth]{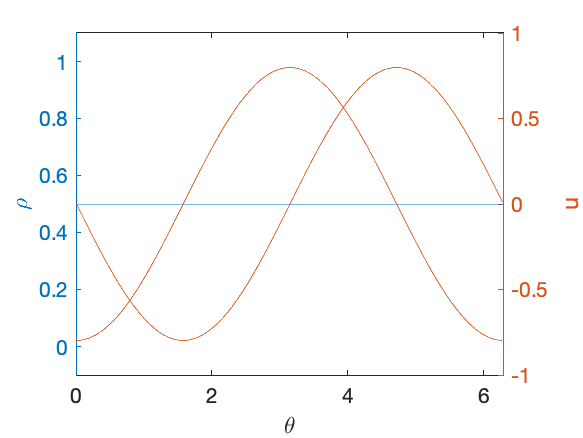}
\includegraphics[width=0.34\textwidth]{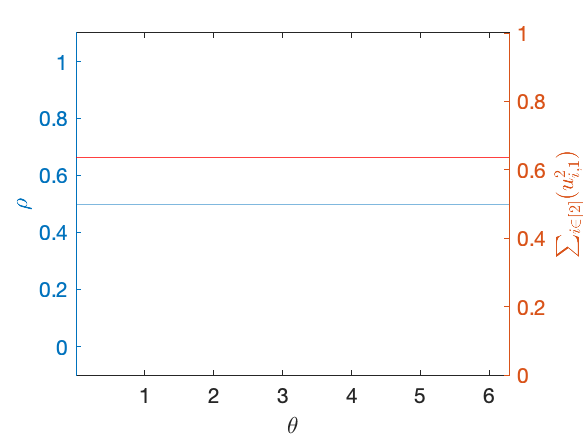}\\
$k=2$ \\
\includegraphics[width=0.34\textwidth]{figs/disc_max_2_2.9193_rho_u.png}
\includegraphics[width=0.34\textwidth]{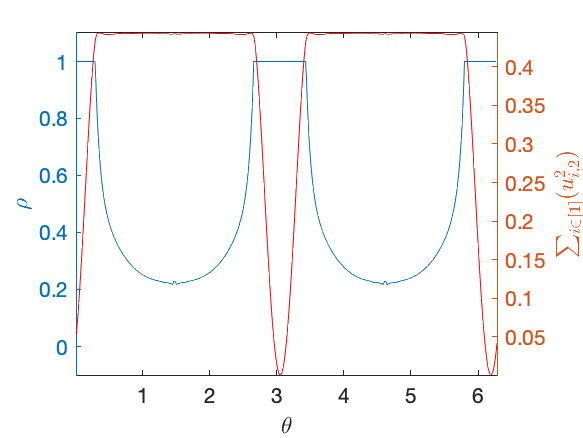}\\
$k=3$ \\
\includegraphics[width=0.34\textwidth]{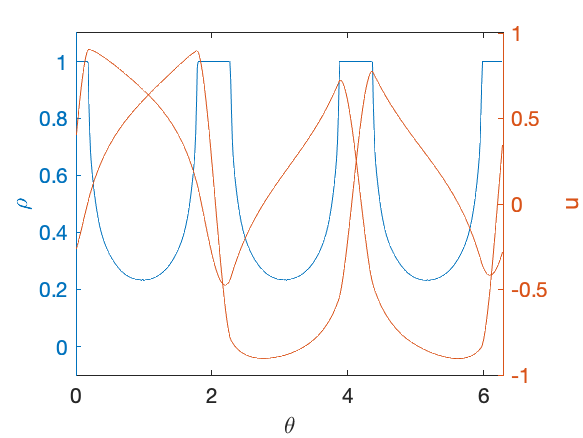}
\includegraphics[width=0.34\textwidth]{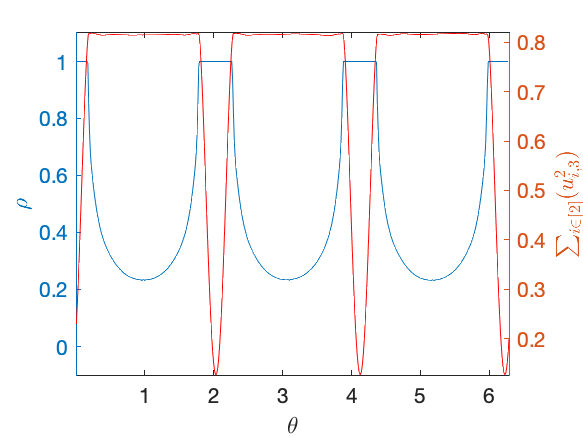}\\
$k=4$ \\
\includegraphics[width=0.34\textwidth]{figs/disc_max_4_5.8386_rho_u.png}
\includegraphics[width=0.34\textwidth]{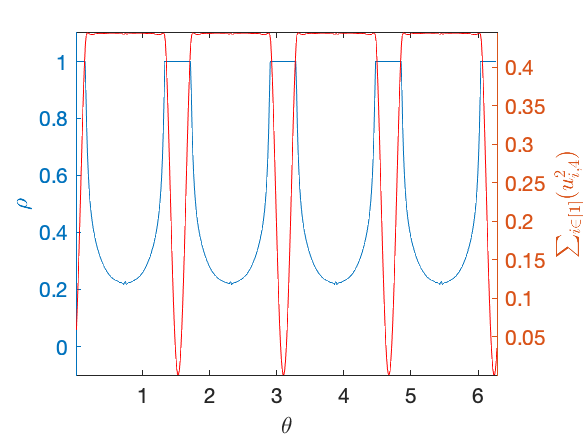}
\caption{Solutions for \eqref{e:RelaxMax} 
with $\Omega = \D$, $\alpha = 0.5$, and $k = 1,2,3,4$ (values indicated). 
In the left hand panels, the maximizing density $\rho_k^\triangle$ is plotted in blue and the corresponding eigenfunction(s) are plotted in red.
For $k$ even, the eigenvalue is simple and the right hand panel shows that the optimality condition \eqref{e:maxRhoCondMult} is satisfied. 
For $k$ odd, the eigenvalue has multiplicity $2$ and the right hand panel shows that the optimality condition \eqref{e:maxRhoCondMult} is satisfied for the two eigenfunctions plotted in the left hand panel. 
See Section~\ref{s:DiskMax}.}
\label{f:SmallkMaximizer}
\end{figure}

For $k\geq 2$, the maximizer $\rho_k^{\triangle}$ in \eqref{e:RelaxMax} has a non-trivial form shown in Figures~\ref{f:EvenkMaximizer} and \ref{f:SmallkMaximizer}. It is the $k$-repetition of an interval that is constant ($=1$) on a first subinterval and is strictly positive and convex on a second subinterval. In the right panels of Figure~\ref{f:EvenkMaximizer}, we identify the Steklov boundary in red. As the maximizers are not bang-bang, the total length of red curves is not necessarily $\pi$ while the constraint $\int_0^{2\pi} \rho(\theta)\, d\theta = \pi$ is satisfied.

For even $k$, $\lambda_{k,\alpha}^{\triangle}$ have multiplicity one. 
Figure~\ref{f:EvenkMaximizer} shows the optimal $\rho_k^\triangle$ and associated eigenfunctions for $\alpha = 0.5$ and $k = 2, 4$. 
In Figure~\ref{f:SmallkMaximizer},  we see that the eigenfunction satisfies the optimality condition \eqref{e:maxRhoCondMult} in Theorem~\ref{prop:OptCondMultEig}. Interestingly, $u_k \Big|_{\Gamma} (\theta)$ is constant on intervals where the optimal density $\rho_k^\triangle$ takes values in $(0,1)$. 
Finally, the optimal values are given by  
\[ \lambda_{k,\alpha}^{\triangle} 
\ = \   \frac{k}{2} \lambda_{2,\alpha}^{\triangle}, 
\qquad \qquad \qquad k\geq 2 \textrm{ even}, \ \alpha \in (0,1]. \]

For odd $k$, $\lambda_{k,\alpha}^{\triangle}$ have multiplicity two. 
Figure~\ref{f:SmallkMaximizer} shows the optimal $\rho_k^\triangle$ and associated eigenfunctions for $\alpha = 0.5$ and $k=3$. 
The corresponding eigenfunctions satisfy the optimality condition \eqref{e:maxRhoCondMult} in Theorem~\ref{prop:OptCondMultEig}. That is, when $\rho_k^\triangle$ takes values in $(0,1)$, for a particular choice of eigenfunctions,  $\sum_{i \in [2]} u^2_{i,k} \Big|_{\Gamma} (\theta) = c$, for some $c>0$. A choice of $u_{i,k}$ is shown in the left panels of Figure~\ref{f:SmallkMaximizer}.  
Finally, the optimal values satisfy 
$ \lambda_{k,\alpha}^\triangle 
 \searrow   
 \frac{k}{2} \lambda_{2,\alpha}^{\triangle}$, 
$k\geq 3$ odd, $\alpha \in (0,1]$. 

In Figure~\ref{f:Weyl}, we plot $\alpha$ vs. $\alpha \lambda_{2,\alpha}^\triangle$. This indicates that, for $k$ even, except in the limit that $\alpha \to 0$, we have that $\lambda_{2,\alpha}^\triangle < \frac{2}{\alpha}$, implying that $\lambda_{k,\alpha}^\triangle = \frac{k}{2} \lambda_{2,\alpha}^\triangle < \frac{k}{\alpha}  =  \frac{2 \pi}{\alpha\abs{\Gamma}} k $, which shows that the Hersch-Payne-Schiffer inequality \cite{Hersch_1974,Girouard_2010} is not strict for the extremal weighted Steklov eigenvalue problem (and hence the extremal S-N eigenvalue problem) for $k$ even and $\alpha > 0$. 

We remark that the oscillations for the maximizing densities, $\rho_k^\triangle$, are reminiscent of the maximizing domains with oscillatory boundaries for the shape optimization problem studied in \cite{Akhmetgaliyev2016}.

\subsection*{Comparison of \texorpdfstring{$\lambda_{k,\alpha}^{\bigtriangledown}$}{λ{k,α}▽} and \texorpdfstring{$\lambda_{k,\alpha}^{\triangle}$}{λ{k,α}△} to Weyl's law for the disk.} 
\begin{figure}[t]
\centering
\includegraphics[width=.5\linewidth]{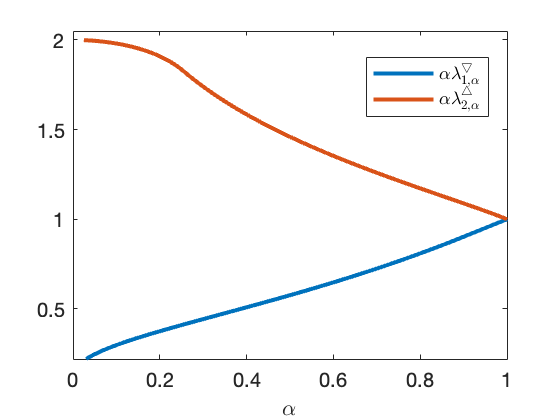}
\caption{A plot of  $\alpha$ vs. $\alpha \lambda_{1,\alpha}^{\bigtriangledown}$ and $\alpha \lambda_{2,\alpha}^{\triangle}$ for a comparison of the extremal eigenvalues with Weyl's law; see \eqref{e:CompWeyl}.}
\label{f:Weyl}
\end{figure}

We recall that Weyl's law on a unit disk gives that for any \emph{fixed} $\rho \in \calA_\alpha$, 
$\lambda_k \sim \frac{2\pi}{\alpha\abs{\Gamma}} \left\lceil \frac{k}{2} \right\rceil  = \frac{1}{\alpha} \left\lceil \frac{k}{2} \right\rceil =:\lambda_{k,\alpha}^W.$
We now compare the values of  $\lambda_{k,\alpha}^{\bigtriangledown}$ and $\lambda_{k,\alpha}^{\triangle}$ (corresponding to sequences $\{\rho_{k,\alpha}^{\bigtriangledown}\}_{k\geq 1}$ and $  \{\rho_{k,\alpha}^{\triangle} \}_{k\geq 1}\subset \calA_\alpha$) to the value ``predicted'' by Weyl's law, $\lambda_{k,\alpha}^W$. From our previous numerical experiments, we have, as $k\to \infty$,
\begin{equation}
\label{e:CompWeyl}
\frac{\lambda_{k,\alpha}^{\bigtriangledown}}{\lambda_{k,\alpha}^W} \to 
 \alpha \sigma_{1,\alpha}^\bigtriangledown = \alpha \sigma_1(\calI_\alpha^{2})
\qquad \qquad \textrm{and} \qquad \qquad 
\frac{\lambda_{k,\alpha}^{\triangle}}{\lambda_{k,\alpha}^W} \to \alpha \lambda_{2,\alpha}^\triangle. 
\end{equation} 
In Figure~\ref{f:Weyl}, we plot $\alpha$ vs. $\alpha\lambda_{1, \alpha}^\bigtriangledown = \alpha \sigma_{1,\alpha}^\bigtriangledown$  and $\alpha \lambda_{2,\alpha}^\triangle$. We see that for $\alpha$ small, say $\alpha =0.02$, we have that $\alpha \sigma_{1,\alpha}^\bigtriangledown \approx 0.2$ which indicates that 
$\sigma_k(\calI_\alpha^{k+1})$ is only $20\%$ what Weyl's law would predict! 
Furthermore, for $\alpha$ small, we have that $\alpha \lambda_{2,\alpha}^\triangle \approx 2$  which indicates that 
$\lambda_{k,\alpha}^{\triangle}$ is twice what Weyl's law would predict for large $k$. 
For other examples of where extremal eigenvalues are a multiplicative factor different than predicted by Weyl's law, see \cite{Girouard_2010,Akhmetgaliyev2016,Kao_2017,Karpukhin_2021,Kao_2023}.

\printbibliography

\clearpage 
\appendix 
\section{A counterexample to convexity properties of higher eigenvalues} \label{a:counterexample}
We consider the eigenvalue problem 
\[ \begin{cases}
    \, -\Delta u(x) = \lambda\rho(x)u(x) & x\in\Omega, \\
    \, u(x) = 0 & x\in\del\Omega,
\end{cases} \]
where $\rho$ is in the admissible set 
\[ \calA\coloneqq \left\{\rho\in L^{\infty}(\Omega)\colon \alpha\le \rho(x)\le \beta \textrm{ for a.e. } x\in\Omega, \, \int_{\Omega} \rho(x)\, dx = c\right\}, \]  
with $0\le \alpha < \beta$. 
In \cite[Theorem 9.1.3]{Henrot_2006}, it is claimed that $\rho\mapsto \lambda^{-1}_k(\rho)$ is convex for fixed $k\geq 1$, but this is not true for $k\geq 2$ (as the author agreed in personal communication).

\begin{prop}
For $k\geq 2$, the map $\rho \mapsto \lambda^{-1}_k(\rho)$ is not convex on $\mathcal{A}$. 
\end{prop}
\begin{proof}
For simplicity, a one-dimensional counterexample is given here. We consider the eigenvalue problem 
\[ \begin{cases} 
\, -u^{\prime\prime}(x) = \lambda\rho(x)u(x), & x \in (-1,1), \\
\, u(-1) = u(1) = 0, 
\end{cases} \]
with
$ \rho(x) = \begin{cases}
\, \rho_-, & \textrm{if } x < 0, \\
\, \rho_+,  & \textrm{if } x \geq 0.
\end{cases} $ 
Note that this is a Sturm-Liouville problem and has discrete spectrum consisting only of simple eigenvalues $\{\lambda_k\}_{k = 1}^\infty$. 
We seek an eigenfunction of the form 
\[ u(x) = \begin{cases}
\, C_- \sin\left(\sqrt{\lambda\rho_-}(1+x)\right), & \textrm{if } x < 0, \\
\, C_+ \sin\left(\sqrt{\lambda\rho_+}(1-x)\right), & \textrm{if } x \geq 0. 
\end{cases} \] 
Imposing continuity of $u$ and $u'$ at $x = 0$ yields the matrix equation $AC = 0$, where 
$A = \begin{pmatrix}
            \sin(\sqrt{\lambda\rho_-}) & -\sin(\sqrt{\lambda\rho_+}) \\ 
            \sqrt{\rho_-} \cos(\sqrt{\lambda\rho_-}) & \sqrt{\rho_+}\cos (\sqrt{\lambda\rho_+}) 
        \end{pmatrix} 
$ 
and 
$C = \begin{pmatrix} C_- \\  C_+ \end{pmatrix}$.
Setting $\det(A) = 0$, we obtain a transcendental equation for $\lambda$, 
$ \sqrt{\rho_+} \tan(\sqrt{\lambda\rho_-}) + \sqrt{\rho_-} \tan(\sqrt{\lambda\rho_+}) = 0$.
For fixed $\rho_-, \rho_+$, the nonzero roots $\lambda_k = \lambda_k(\rho_-, \rho_+)$, $k \in \mathbb N$ are the eigenvalues. 

We now consider  
$ \rho_{1}(x) = 
    \begin{cases}
        \, 1, & \textrm{if } x < 0, \\
        \, 4, & \textrm{if } x\ge 0,
        \end{cases}$ 
and
$\rho_{2}(x) = 
    \begin{cases}
    \, 4, & \textrm{if } x < 0, \\
    \, 1, & \textrm{if } x\ge 0, 
    \end{cases}$
and let $\rho$ denote the convex combination of $\rho_{1}$ and $\rho_{2}$, i.e., $\rho(x; t) = t\rho_{1}(x) + (1 - t)\rho_{2}(x)$ with convex parameter $t\in[0,1]$. In particular, we have 
\begin{align*}
\rho_-(x; t) & = t + 4(1 - t) = 4 - 3t, \\  
\rho_+(x; t) & = 4t + (1 - t) = 1 + 3t. 
\end{align*}
Moreover, $\int_{-1}^1 \rho(x; t)\, dx = 5$. Next we define the function 
\[ f(\lambda; t)\coloneqq \sqrt{1 + 3t}\tan(\sqrt{\lambda(4 - 3t)}) + \sqrt{4 - 3t}\tan(\sqrt{\lambda (1 + 3t)}). \] 
The eigenvalues $\lambda_k(t)$ satisfy $f(\lambda(t); t) = 0$. For $t = 0$, the eigenvalues satisfy 
\begin{equation} \label{e:SLeigs} 
2\tan(\sqrt{\lambda}) = - \tan\left(2\sqrt{\lambda}\right) = \frac{-2\tan(\sqrt{\lambda})}{1 - \tan^2(\sqrt{\lambda})}, 
\end{equation} 
which gives that either $\tan(\sqrt{\lambda}) = 0$ or $\tan(\sqrt{\lambda}) = \pm \sqrt{2}$. 
One can verify that for $k\ge 1$, 
\[ \lambda_k(t = 0) = \begin{cases}
\left(\left\lfloor \frac{k}{3}\right\rfloor \pi+\tan^{-1}\left(\sqrt{2}\right)\right)^{2}, & \textrm{if} \mod(k,3) = 1,\\
\left(\left\lfloor \frac{k}{3}\right\rfloor \pi+\pi-\tan^{-1}\left(\sqrt{2}\right)\right)^{2}, & \textrm{if} \mod(k, 3) = 2,\\
\left(\left\lfloor \frac{k}{3}\right\rfloor \pi\right)^{2}, & \textrm{if} \mod(k, 3) = 0. 
\end{cases}
\]

We will show that the map $t\mapsto \lambda^{-1}_{2} (t)$ is not convex on $[0,1]$. Note that $f(\lambda,1-t) = f(\lambda,t)$ which implies that $\lambda_{k}(1-t) = \lambda_{k}(t)$, which is just invariance of the eigenvalues to the reflection $x\mapsto -x$.  
We show that $\dot\lambda_2(t=0) < 0$. We invoke the implicit function theorem and use the formula 
\[ \dot\lambda_k(t = 0) = \frac{d\lambda_k}{dt}\bigg|_{t=0} = - \frac{f_t}{f_\lambda}\bigg|_{t=0,\, \lambda = \lambda_k(0)}. \] 
For notational simplicity, let us write $\eta_k = \lambda_k(0)$. Differentiating with respect to $\lambda$, we obtain 
\begin{align*}
f_\lambda\big|_{t = 0,\, \lambda = \eta_k} = & \frac{1}{\sqrt{\eta_k}}\left[\sec^2(\sqrt{4\eta_k}) + \sec^2(\sqrt{\eta_k})\right] > 0. 
\end{align*} 
Differentiating with respect to $t$, we obtain 
\begin{equation}\label{e:SL1} 
f_t\big|_{t = 0,\, \lambda = \eta_k} = -\frac{3}{4}\left(- 2\tan(2\sqrt{\eta_k}) + \sqrt{\eta_k}\sec^2(2\sqrt{\eta_k}) + \tan(\sqrt{\eta_k}) - 4\sqrt{\eta_k}\sec^2(\sqrt{\eta_k})\right). 
\end{equation} 
To simplify \eqref{e:SL1}, we square both sides of \eqref{e:SLeigs} with $\lambda = \eta_k$ and rewrite using the identity $\tan^2(z) = \sec^2(z) - 1$: 
\[ 4\tan^2(\sqrt{\eta_k}) = \tan^2(2\sqrt{\eta_k}) \implies 4\sec^2(\sqrt{\eta_k}) = \sec^2(2\sqrt{\eta_k}) + 3. \]
Consequently, \eqref{e:SL1} simplifies to 
\[ f_t\big|_{t = 0, \, \lambda = \eta_k} = -\frac{3}{4}\left(-3\sqrt{\eta_k} + 5\tan(\sqrt{\eta_k})\right). \]
We now take $k = 2$ so that $\eta_2 = \lambda_2(0) = \left(\pi - \tan^{-1}(\sqrt{2})\right)^2$ and $\tan(\sqrt{\lambda_2(0)}) = - \sqrt{2}$. We then have that $f_t\big|_{t = 0,\, \lambda = \lambda_2(0)}>0$, which gives $\dot\lambda_{2}(t = 0) < 0$ and $\frac{d}{dt}(\lambda^{-1}_{2})\big|_{t = 0} > 0$. This contradicts the convexity of $\lambda^{-1}_{2}$. 
\end{proof}

Finally, we remark that convexity/concavity properties for higher eigenvalues of other operators is also misstated, e.g., 
\cite[Theorem 8.1.2]{Henrot_2006} for the Schr\"odinger equation and 
\cite[Theorem 10.1.1]{Henrot_2006} for the conductivity equation.

\end{document}